\documentclass[10pt,a4paper,twoside,openany]{report}
\usepackage[cp1251]{inputenc} 
\usepackage[ukrainian]{babel} 
\usepackage{amsfonts,amsmath,amssymb,amscd}
\usepackage{bezier,graphics,graphicx}
\usepackage{epsfig, layout,latexsym}
\usepackage{theorem}
\usepackage{cite}
\usepackage{indentfirst}
\usepackage{array}

\AtBeginDocument{} \AtBeginDocument{}

\theorembodyfont{\sl}
\newtheorem{theorem}{Theorem}[section]

\newtheorem{lemma}[theorem]{Lemma}
\newtheorem{corollary}[theorem]{Corollary}
\newtheorem{remark}[theorem]{Remark}

\numberwithin{equation}{section} \numberwithin{theorem}{section}
\newcommand{\card}{\text{card}\medskip}

\makeatletter
\renewenvironment{thebibliography}[1]
 {\section*{\centerline{\rm\textsc{Bibliography}}}%
 \@mkboth{\MakeUppercase\refname}{\MakeUppercase\refname}%
 \list{\@biblabel{\@arabic\c@enumiv}}%
 {\settowidth\labelwidth{\@biblabel{#1}}%
 \leftmargin\labelwidth
 \advance\leftmargin\labelsep
 \@openbib@code
  \usecounter{enumiv}%
  \let\p@enumiv\@empty
  \renewcommand\theenumiv{\@arabic\c@enumiv}}%
  \sloppy
  \clubpenalty4000
  \@clubpenalty \clubpenalty
  \widowpenalty4000%
  \sfcode`\.\@m
  \setlength{\itemsep}{-0.1cm}}
  {\def\@noitemerr
  {\@latex@warning{Empty 'thebibliography' environment}}%
 \endlist}
\renewcommand{\@biblabel}[1]{#1.}

\makeatother \hfuzz=0.5pt \tolerance=500

\begin{document}

\setcounter{equation}{0} \setcounter{figure}{0} \setcounter{table}{0} \setcounter{footnote}{0} \setcounter{section}{0}

\begin{center}
\textbf{On optimal recovering high order partial derivatives of bivariate functions }

\end{center}

\def\headlinetitle{Error bounds for Fourier-Legendre truncation method}

\vspace{-4mm} \noindent


\vspace*{5mm}
\centerline{\textsc { Y.V. Semenova $\!\!{}^{\dag}$, S.G. Solodky  $\!\!{}^{\dag,\ddag}$} }


%
 \centerline{$\!\!{}^{\dag}\!\!$ Institute of Mathematics, National Academy of Sciences of Ukraine, Kyiv}
\centerline{$\!\!{}^{\ddag}\!\!$ University of Giessen, Department of Mathematics, Giessen, Germany}

%
%

\vspace{4mm}
\begin{small}
\begin{quote}
\textsc{Abstract.}
The problem of recovering partial derivatives of high orders of bivariate functions with finite
smoothness is studied. Based on the truncation method, a numerical differentiation algorithm was constructed, 
which is optimal by the order both in the sense of accuracy and in the sense of the amount of Galerkin information involved. 
Numerical demonstrations are provided to illustrate that the proposed method can be implemented successfully.
\end{quote}
\end{small}

\begin{small}
	\begin{quote}
	\end{quote}
\end{small}


\section{Introduction}

The present work is devoted to the issue of research and improving the efficiency of numerical differentiation methods.
As is known, the study of numerical differentiation from the theory of ill-posed problems standpoint originates from the work \cite{Dolgopolova&Ivanov_USSR_Comput_Math_Math_Phys_1966_Eng} and by now there are many different methods to the stable recovery of derivatives with respect to perturbed data (see, for example, \cite{Ramm_1968_No11}, \cite{VasinVV_1969_V7_N2},
\cite{Groetsch_1992_V74_N2}, \cite{Hanke&Scherzer_2001_V108_N6}, \cite{Ahn&Choi&Ramm_2006}, \cite{Lu&Naum&Per},
\cite{Nakamura&Wang&Wang_2008}, \cite{Zhao_2010}, \cite{Zhao&Meng&Zhao&You&Xie_2016}, \cite{Meng&Zhaoa&Mei&Zhou_2020},
\cite{SSS_CMAM}). As to provide stability of approximation, all mentioned papers can be fell into two directions.

The standard approach to ensuring the stability of a numerical differentiation problem is to apply the classical Tikhonov regularization methods with an appropriate selection of the regularization parameter (see, in particular, \cite{Dolgopolova&Ivanov_USSR_Comput_Math_Math_Phys_1966_Eng}, \cite{VasinVV_1969_V7_N2}, \cite{Cul71}, \cite{Zhao&Meng&Zhao&You&Xie_2016}). The main efforts of researchers in this direction are focused on determining the optimal value of the regularization and discretization parameters, which is often a non-trivial task.

However, in some cases, there is an alternative approach to achieving stability. This approach is called self-regularization and consists in choosing the appropriate discretization level depending on the noise level of the input data. Here, the discretization level acts as a regularization parameter, due to which stability is ensured. Examples of using self-regularization to solve some classes of ill-posed problems can be found in \cite{VainHam}, \cite{SolSem2012}, \cite{HamKan2018}. As for the numerical differentiation problem, the idea of using self-regularization to develop stable algorithms was earlier proposed in \cite{Groetsch_1991} and \cite{RammSmir_2001}. In this paper, we will continue to study the approximation properties of self-regularization and focus on recovering the partial derivatives of bivariate functions by finite Fourier sums (spectral truncation method) and provide a stable approximation with an appropriate choice of summation limit. The method was first applied to the problem of numerical differentiation in \cite{Lu&Naum&Per}. Subsequently, the effectiveness of this approach was confirmed by the results of \cite{Sem_Sol_2021}, \cite{Sol_Stas_UMZ2022}, \cite{SSS_CMAM}.
Continuing the series of studies devoted to the spectral truncation method, in this paper, for the approximation of partial derivatives, a modification of the spectral truncation method will be presented in combination with a priori choice of the discretization level depending on the noise. At the same time, the approach under consideration is expected not only to ensure the optimal accuracy of approximations but also the efficiency of the usage of computing resources.

The article is organized as follows. In Section 2, the necessary definitions are introduced and the problem statement for optimizing numerical differentiation methods in the sense of the minimal Galerkin information radius is given. In Sections 3 and 4, for 
a proposed version of the spectral truncation method its error estimates in quadratic and uniform metrics, respectively, are established. Section 5 is devoted to finding order estimates for the minimal radius of Galerkin information while establishing the optimality (on a power scale) of the spectral truncation method studied above. Section 6 presents the results of numerical experiments demonstrating the effectiveness of the method under consideration.

\section{Description of the problem}

Let $\{\varphi_k(t)\}_{k=0}^\infty$  be the system of Legendre polynomials  orthonormal on $[-1,1]$ as
$$
\varphi_k(t)=\sqrt{k+1/2}(2^kk!)^{-1}\frac{d^k}{dt^k}[(t^2-1)^k], \quad k=0,1,2,\ldots .
 $$

 By $L_2=L_2(Q)$   we mean space of square-summable  on $Q=[-1,1]^2$ functions $f(t,\tau)$ with inner  product
$$
\langle f, g\rangle=\int_{-1}^{1}\int_{-1}^{1}f(t,\tau)g(t,\tau)d \tau d t
$$
and standard norm
$$
\|f\|_{L_2}^2=\sum_{k,j=0}^{\infty}|\langle f, \varphi_{k,j} \rangle|^2 < \infty ,
$$
where $$ \langle f, \varphi_{k,j}\rangle=\int_{-1}^{1}\int_{-1}^{1}f(t,\tau)\varphi_k(t)\varphi_j(\tau)d\tau dt, \quad
k,j=0,1,2,\ldots,
$$
are Fourier-Legendre coefficients of $f$. Let $\ell_p$, $1\leq p\leq\infty$, be the space of numerical sequences
$\overline{x}=\{x_{k,j}\}_{k,j\in\mathbb{N}_0}$, $\mathbb{N}_0=\{0\}\bigcup\mathbb{N}$, such that the corresponding
relation
$$
\|\overline{x}\|_{\ell_p}  := \left\{
\begin{array}{cl}
\bigg(\sum\limits_{k,j\in\mathbb{N}_0} |x_{k,j}|^p\bigg)^{\frac{1}{p}} < \infty ,
 \ & 1\leq p<\infty ,
\\\\
\sup\limits_{k,j\in\mathbb{N}_0}  |x_{k,j}| < \infty ,
  \ & p=\infty ,
\end{array}
\right.
$$
is fulfilled.

We introduce   the space of functions
$$
 L_{s,2}^{\overline{\mu}} := L_{s,2}^{\overline{\mu}}(Q) := \{f\in L_2(Q)\!: \|f\|_{s,\overline{\mu}}^{s} =
 \sum\limits_{k,j=0}^{\infty} {\underline{k}}^{s\mu_1} {\underline{j}}^{s\mu_2}
 |\langle f , \varphi_{k,j} \rangle |^s < \infty
   \} ,
$$
where $\overline{\mu}=(\mu_1,\mu_2)$, $\mu_1, \mu_2 > 0$,\ $1\le s<\infty$,\ $\underline{k}=\max\{1,k\}$,\
$k=0,1,2,\dots$. Note that in the future we will use the same notations both for space and for a unite ball from this
space: $L_{s,2}^{\overline{\mu}} = L_{s,2}^{\overline{\mu}}(Q) = \{f\in L_{s,2}^{\overline{\mu}}\!:
\|f\|_{s,\overline{\mu}} \leq 1\}$, what we call  a class of functions. What exactly is meant by $
L_{s,2}^{\overline{\mu}}$, space or class will be clear depending on the context in each case. It should be noted that
$L^{\overline{\mu}}_{s,2}$ is a generalization of the class of bivariate functions with dominating mixed derivatives.
Moreover, let $C=C(Q)$ be the space of continuous on $Q$ bivariate functions.

We represent a function $f(t,\tau)$ from $L_{s,2}^{\overline{\mu}}$ as
$$
 f(t,\tau) = \sum_{k,j=0}^{\infty} \langle f, \varphi_{k,j}\rangle \varphi_k(t)\varphi_j(\tau),
$$
and by $r$-th partial derivative of $f$ we mean the following series
\begin{equation}\label{r_deriv}
f^{(r,0)}(t,\tau) =  \sum_{k=r}^{\infty} \sum_{j=0}^{\infty} \langle f, \varphi_{k,j}\rangle
\varphi^{(r)}_k(t)\varphi_j(\tau), \quad r=1,2,\ldots.
\end{equation}
Assume that instead of the exact values of  the Fourier-Legendre coefficients $\langle f, \varphi_{k,j} \rangle$ only
some of their perturbations are known with the error level $\delta$ in the metrics of $\ell_p$, $1\leq p\leq\infty$.
 More accurately, we assume that there is a sequence of numbers $\overline{f^\delta}= \{\langle
f^\delta, \varphi_{k,j} \rangle\}_{k,j\in\mathbb{N}_0}$ such that for $\overline{\xi}=
\{\xi_{k,j}\}_{k,j\in\mathbb{N}_0}$, where $\xi_{k,j}=\langle f-f_\delta,\varphi_{k,j}\rangle$, and for some $1\leq
p\leq \infty$ the relation
 \begin{equation}\label{perturbation}
   \|\overline{\xi}\|_{\ell_p} \leq \delta , \quad 0<\delta <1 ,
 \end{equation}
is true.

The research of this work is devoted to the optimization of methods for recovering the derivative (\ref{r_deriv}) of
functions from the class $L_{s,2}^{\overline{\mu}}$. Further, we give a strict statement of the problem to be studied.
In the coordinate plane $[r,\infty)\times[0,\infty)$ we take an arbitrary bounded domain $\Omega$. By $\card(\Omega)$
we mean the number of points that make up $\Omega$ and by the information vector $G(\Omega, \overline{f}^{\delta})\in
\mathbb{R}^{N}$, $\card (\Omega) = N$ we take the set of perturbed values of Fourier-Legendre coefficients
$\left\{\langle f^{\delta}, \varphi_{k,j} \rangle\right\}_{(k,j)\in\Omega}$.

Let  $X=L_{2}(Q)$ or $X=C(Q)$.
 By numerical differentiation algorithm we mean any mapping $\psi^{(r,0)} = \psi^{(r,0)}(\Omega)$ that corresponds to the information vector
 $G(\Omega, \overline{f}^{\delta})$ an element $\psi^{(r,0)}(G(\Omega, \overline{f}^{\delta})) \in X$,
 which is taken as an approximation to the derivative (\ref{r_deriv}) of function $f$ from $L_{s,2}^{\overline{\mu}}$.
We denote by  $\Psi(\Omega)$ the set of all algorithms $\psi^{(r,0)}(\Omega):\,\mathbb{R}^{N}\rightarrow X$, that use
the same information vector $G(\Omega, \overline{f}^{\delta})$.

In general, algorithms in $\Psi(\Omega)$ are not required to be linear or stable.
The only condition for algorithms in $\Psi(\Omega)$ is to use input information
in the form of perturbed values of the Fourier-Legendre coefficients with indices from the domain $\Omega$  of the
coordinate plane. Such a general understanding of the algorithm is explained by the desire to compare the widest
possible range of possible methods for numerical differentiation.

The error of the algorithm $\psi^{(r,0)}$ on the class $L^{\overline{\mu}}_{s,2}$ is determined by the quantity
 $$
 \varepsilon_{\delta}(L^{\overline{\mu}}_{s,2}, \psi^{(r,0)}(\Omega), X, \ell_p)
 = \sup\limits_{f\in L^{\overline{\mu}}_{s,2}, \, \|f\|_{s,\overline{\mu}}\leq 1}
 \ \, \sup\limits_{\overline{f}^{\delta}: \, (\ref{perturbation})}
 \| f^{(r,0)} - \psi^{(r,0)}(G(\Omega, \overline{f}^{\delta})) \|_X .
 $$

The minimal radius of the Galerkin information for the problem of numerical differentiation on the class
$L^{\overline{\mu}}_{s,2}$ is given by
$$
R^{(r,0)}_{N,\delta} (L^{\overline{\mu}}_{s,2}, X, \ell_p) = \inf\limits_{\Omega: \, \card(\Omega)\leq N}
 \ \, \inf\limits_{\psi^{(r,0)}\in\Psi(\Omega)} \varepsilon_{\delta}(L^{\overline{\mu}}_{s,2}, \psi^{(r,0)}(\Omega), X, \ell_p) .
$$

The value  $R^{(r,0)}_{N,\delta} (L^{\overline{\mu}}_{s,2}, X, \ell_p)$ describes the minimal possible accuracy in the
metric of space $X$, which can be achieved by numerical differentiation of arbitrary function $f\in
L_{s,2}^{\overline{\mu}}$ , while using not more than $N$ values of its Fourier-Legendre coefficients that are
$\delta$-perturbed in the $\ell_p$ metric. Note that the minimal radius of Galerkin information in the problem of
recovering the first partial derivative was studied in \cite{Sol_Stas_UMZ2022}, and for other types of ill-posed
problems, similar studies were previously carried out in \cite{PS1996}, \cite{Mileiko_Solodkii_2014}. It should be
added that the minimal radius characterizes the information complexity of the considered problem and is traditionally
studied within the framework of the IBC-theory (Information-Based Complexity Theory), the foundations of which are laid
in monographs \cite{TrW} and \cite{TrWW}.

Our research aim is to construct a numerical differentiation method that achieves sharp estimates
(in the power scale) for the quantities $R^{(r,0)}_{N,\delta}
(L^{\overline{\mu}}_{s,2}, C, \ell_p)$ and $R^{(r,0)}_{N,\delta} (L^{\overline{\mu}}_{s,2}, L_{2}, \ell_p)$.

\section{Truncation method. Error estimation in $L_2$ metric }

As an algorithm for the numerical differentiation of functions from $L^{\overline{\mu}}_{s,2}(Q)$ we will use
a truncation method. The essence of this method is to replace the Fourier series (\ref{r_deriv}) with a finite
Fourier sum using perturbed data $\langle f^\delta, \varphi_{k,j} \rangle$. In the truncation method to ensure the
stability of the approximation and achieve the required order of accuracy,  it is necessary to choose properly the
discretization parameter, which here serves as a regularization parameter. So, the process of regularization in the method
under consideration consists in matching the discretization parameter with the perturbation level $\delta$ of the input
data. The simplicity of implementation is the main advantage of this method.

In the case of an arbitrary bounded domain $\Omega$ of the coordinate plane $[r,\infty)\times[0,\infty)$, the
truncation method for differentiation of bivariate functions has the form
$$
\mathcal{D}_\Omega^{(r,0)} f^\delta(t,\tau)=\sum_{(k,j)\in \Omega} \langle f^\delta, \varphi_{k,j} \rangle
\varphi_k^{(r)}(t)\varphi_j(\tau).
$$
To increase the efficiency of the approach under study, as the domain $\Omega$ we take a hyperbolic cross of
the following form
$$
\Omega = \Gamma_{n,\gamma} := \{ (k,j):\ k\cdot j^{\gamma} \leq n,\quad k=r,\ldots, n,\quad j=0,\ldots,
(n/r)^{1/\gamma}\}, \qquad \gamma\geq 1 .
$$
Then the proposed version of the truncation method can be written as
\begin{equation} \label{ModVer}
\mathcal{D}_{n,\gamma}^{(r,0)} f^\delta(t,\tau) = \sum_{k=r}^{n} \sum_{j=0}^{(n/k)^{1/\gamma}} \langle f^\delta,
\varphi_{k,j}\rangle \varphi^{(r)}_k(t)\varphi_j(\tau).
\end{equation}
We note that earlier the idea of a hyperbolic cross in the problem of numerical differentiation was used in the papers
\cite{Sem_Sol_2021}, \cite{Sol_Stas_UMZ2022}, \cite{SSS_CMAM} (for more details about usage of hyperbolic cross in
solving the other types of ill-posed problems see \cite{Pereverzev_Computing_1995}, \cite{Sem_Sol_2008},
\cite{Mileiko_Solodkii_2017_UMJ}, \cite{Sol_Stas_JC2020}).

The approximation properties of the method (\ref{ModVer}) will be investigated in Sections 3 and 4 while  in Section 5
it will be found that the method (\ref{ModVer}) is order-optimal in the sense of the minimal radius of the Galerkin
information.

The parameters $n$ and $\gamma$ in (\ref{ModVer}) should be chosen depending on $\delta$ and $\overline{\mu}$ so as to minimize
 the error of the method $\mathcal{D}_{n,\gamma}^{(r,0)}$:
\begin{equation}\label{fullError}
f^{(r,0)}(t,\tau)-\mathcal{D}_{n,\gamma}^{(r,0)} f^\delta(t,\tau)=
\left(f^{(r,0)}(t,\tau)-\mathcal{D}_{n,\gamma}^{(r,0)} f(t,\tau)\right)+\left(\mathcal{D}_{n,\gamma}^{(r,0)}
f(t,\tau)-\mathcal{D}_{n,\gamma}^{(r,0)} f^\delta(t,\tau)\right).
\end{equation}
For the first difference on the right-hand side of (\ref{fullError}), the representation holds true
\begin{equation}\label{Bound_err}
 f^{(r,0)}(t,\tau)-\mathcal{D}_{n,\gamma}^{(r,0)} f(t,\tau)= \triangle_{1}(t,\tau)+\triangle_{2}(t,\tau) ,
\end{equation}
where
\begin{equation}\label{Triangle_1HC}
\triangle_{1}(t,\tau)= \sum_{k=n+1}^{\infty} \, \sum_{j=0}^{\infty} \langle f,
\varphi_{k,j}\rangle\varphi^{(r)}_k(t)\varphi_j(\tau),
\end{equation}
\begin{equation}\label{Triangle_2HC}
 \triangle_{2}(t,\tau)= \sum_{k=r}^{n} \, \sum_{j>(n/k)^{1/\gamma}} \langle f,
 \varphi_{k,j}\rangle\varphi^{(r)}_k(t)\varphi_j(\tau) .
\end{equation}


For our calculations, we need the following formula (see Lemma 18 \cite{Mul69})
\begin{equation}\label{Muller}
 \varphi_{k}'(t) \, = 2 \, \sqrt{k+1/2}
  \mathop{{\sum}^*}\limits_{l=0}^{k-1} \sqrt{l+1/2} \, \varphi_{l}(t) ,
  \ \ \ k\in\mathbb{N} ,
 \end{equation}
where in aggregate \quad $ \mathop{{\sum}^*}\limits_{l=0}^{k-1} \sqrt{l+1/2} \, \varphi_{l}(t) \, $ \quad the summation
is extended over only those terms for which $k+l$ is odd.

Let us estimate the error of the method (\ref{ModVer}) in the metric of $L_2$. A bound for difference (\ref{Bound_err})
is contained in the following statement.

\begin{lemma}\label{lemma_BoundErrHC}
Let $f\in L^{\overline{\mu}}_{s,2}$, $2\leq s< \infty$, $\mu_1>2r-1/s+1/2$, $\mu_2>\mu_1-2r$.
 Then for $1\leq \gamma < \frac{\mu_2+1/s-1/2}{\mu_1-2r+1/s-1/2}$ it holds
$$
  \|f^{(r,0)}-\mathcal{D}_{n,\gamma}^{(r,0)} f\|_{L_2}\leq c\|f\|_{s,\overline{\mu}} n^{-\mu_1+2r-1/s+1/2} .
$$
\end{lemma}

\textit{Proof.} Using the formula (\ref{Muller}), from (\ref{Triangle_1HC}) we have
$$
\triangle_{1}(t,\tau)= 2 \Big(\sum_{k=n+1}^{\infty} \sum_{j=0}^{\infty} \sqrt{k+1/2} \:\langle f, \varphi_{k,j}\:
\rangle \, \varphi_{j}(\tau) \mathop{{\sum}^*}\limits_{l_1=r-1}^{k-1}  \sqrt{l_1+1/2} \,
\varphi'_{l_1}(t)\Big)^{(r-2)}_t
$$
$$
= 4 \Big(\sum_{k=n+1}^{\infty} \sum_{j=0}^{\infty} \sqrt{k+1/2} \:\langle f, \varphi_{k,j}\: \rangle \,
\varphi_{j}(\tau) \mathop{{\sum}^*}\limits_{l_1=r-1}^{k-1}  (l_1+1/2)\: \mathop{{\sum}^*}\limits_{l_2=r-2}^{l_1-1}
\sqrt{l_2+1/2} \, \varphi'_{l_2}(t)\Big)^{(r-3)}_t
$$
$$
= \ldots = 2^r \sum_{k=n+1}^{\infty} \sum_{j=0}^{\infty} \sqrt{k+1/2} \:\langle f, \varphi_{k,j}\:\rangle
 \, \varphi_{j}(\tau)
$$
$$
\times \mathop{{\sum}^*}\limits_{l_1=r-1}^{k-1}  (l_1+1/2)\: \mathop{{\sum}^*}\limits_{l_2=r-2}^{l_1-1} (l_2+1/2)\ldots
\mathop{{\sum}^*}\limits_{l_{r-1}=1}^{l_{r-2}-1} (l_{r-1}+1/2) \mathop{{\sum}^*}\limits_{l_{r}=0}^{l_{r-1}-1}
\sqrt{l_{r}+1/2} \,\varphi_{l_r}(t) .
$$
We note that  in the representation $\triangle_{1}$ only those terms  take place for which all  indexes
$l_1+k, l_2+l_1,...,l_r+l_{r-1}$ are odd. Such a rule is valid also for other terms, namely $\triangle_{2}$
 and $\mathcal{D}_{n,\gamma}^{(r,0)} f-\mathcal{D}_{n,\gamma}^{(r,0)} f^\delta$, appearing in the error representation  (see (\ref{Bound_err}) -- (\ref{Triangle_2HC})).
In the following, for simplicity, we will omit the symbol "*"\, when denoting such summation operations, while taking into account this rule in the calculations.

Further, we change the order of summation and get
$$
\triangle_{1}(t,\tau)= \triangle_{11}(t,\tau)+\triangle_{12}(t,\tau),
$$
where
\begin{equation}\label{triangle_{11}}
\triangle_{11}(t,\tau)= 2^r \sum_{j=0}^{\infty} \varphi_{j}(\tau) \mathop{{\sum}}\limits_{l_r=0}^{n-r+1}
\sqrt{l_r+1/2} \: \varphi_{l_r}(t)
 \sum_{k=n+1}^{\infty} \sqrt{k+1/2} \:\langle f, \varphi_{k,j}\:
 \rangle B^r_{k} ,
\end{equation}
\begin{equation}\label{triangle_{12}}
\triangle_{12}(t,\tau)= 2^r \sum_{j=0}^{\infty} \varphi_{j}(\tau) \mathop{{\sum}}\limits_{l_r=n-r+2}^{\infty}
\sqrt{l_r+1/2}\: \varphi_{l_r}(t) \, \sum_{k=l_r+r}^{\infty}\, \sqrt{k+1/2} \:\langle f, \varphi_{k,j}\:
 \rangle B^r_{k}
\end{equation}
and
\begin{equation}\label{axular1}
B^r_{k}:=\mathop{{\sum}}\limits_{l_1=l_r+r-1}^{k-1}  (l_1+1/2) \mathop{{\sum}}\limits_{l_2=l_r+r-2}^{l_1-1}
(l_2+1/2) \ldots \mathop{{\sum}}\limits_{l_{r-1}=l_r+1}^{l_{r-2}-1}  (l_{r-1}+1/2) \leq c k^{2(r-1)} .
\end{equation}

At first, we consider the case $2<s<\infty.$

Let's bound $\triangle_{11}(t,\tau)$ in $L_2$-norm. Using H\"{o}lder inequality  and (\ref{axular1}),  for
$\mu_1>2r-1/s+1/2$ and $\mu_2>\mu_1-2r$ we have
$$
\|\triangle_{11}\|_{L_2}^2 \leq c \sum_{j=0}^{\infty} \underline{j}^{-2\mu_2} \mathop{{\sum}}\limits_{l_r=0}^{n-r+1}
(l_r+1/2)
$$
$$
\times \left( \sum_{k=n+1}^{\infty} k^{s\mu_1} \underline{j}^{s\mu_2} |\langle f, \varphi_{kj}\rangle |^s\right)^{2/s}
\left( \sum_{k=n+1}^{\infty} k^{-(\mu_1-2r+3/2)s/(s-1)} \right)^{2(s-1)/s}
$$
$$
\leq c  \, n^{-2(\mu_1-2r+1/s-1/2)} \sum_{j=0}^{\infty} \underline{j}^{-2\mu_2} \left( \sum_{k=n+1}^{\infty} k^{s\mu_1}
\underline{j}^{s\mu_2} |\langle f, \varphi_{k,j}\rangle |^s\right)^{2/s}
$$
$$
\leq c \|f\|_{s,\overline{\mu}}^2 \, n^{-2(\mu_1-2r+1/s-1/2)} .
$$
Thus, we find
 $$
 \|\triangle_{11}\|_{L_2}\leq c \|f\|_{s,\overline{\mu}} \, n^{-\mu_1+2r-1/s+1/2} .
 $$
Applying the estimating technique above we can bound the norm of $\triangle_{12}(t,\tau):$
$$
\|\triangle_{12}\|_{L_2}^2 \leq c \, \sum_{j=0}^{\infty} \underline{j}^{-2\mu_2} \left( \sum_{k=r}^{\infty} k^{s\mu_1}
\underline{j}^{s\mu_2} |\langle f, \varphi_{k,j}\rangle |^s\right)^{2/s}
 \mathop{{\sum}}\limits_{l_r=n-r+2}^{\infty} (l_r+1/2)^{-2(\mu_1-2r+3/2)+\frac{2(s-1)}{s}+1}
$$
$$
\leq c \|f\|_{s,\overline{\mu}}^2 \, n^{-2(\mu_1-2r+1/s-1/2)} \Big(\sum_{j=0}^{\infty}
\underline{j}^{-\frac{2s}{s-2}\mu_2}\Big)^{(s-2)/s} \leq c \|f\|_{s,\overline{\mu}}^2 \, n^{-2(\mu_1-2r+1/s-1/2)}.
$$
Combining the estimates for $\triangle_{11}(t,\tau)$ and $\triangle_{12}(t,\tau)$ we obtain
$$
\|\triangle_{1}\|_{L_2}\leq\|\triangle_{11}\|_{L_2}+\|\triangle_{12}\|_{L_2}\leq c \|f\|_{s,\overline{\mu}} \,
n^{-\mu_1+2r-1/s+1/2} .
$$

Now using the formula (\ref{Muller}), from (\ref{Triangle_2HC}) we have
$$
\triangle_{2}(t,\tau) = 2^r \sum_{k=r}^{n}\ \sum_{j>(n/k)^{1/\gamma}} \sqrt{k+1/2} \:\langle f, \varphi_{k,j}\rangle
\varphi_{j}(\tau)
$$
$$
 \times \mathop{{\sum}}\limits_{l_1=r-1}^{k-1}  (l_1+1/2) \mathop{{\sum}}\limits_{l_2=r-2}^{l_1-1}  (l_2+1/2) \ldots
\mathop{{\sum}}\limits_{l_{r-1}=1}^{l_{r-2}-1}  (l_{r-1}+1/2) \mathop{{\sum}}\limits_{l_{r}=0}^{l_{r-1}-1}
\sqrt{l_r+1/2}\: \varphi_{l_r}(t) .
$$
Further, we change the order of summation and get
$$
\triangle_{2}(t,\tau)= \triangle_{21}(t,\tau)+\triangle_{22}(t,\tau)+\triangle_{23}(t,\tau) ,
$$
where
$$
\triangle_{21}(t,\tau) = 2^r \sum_{j=2}^{(n/r)^{1/\gamma}+1} \varphi_{j}(\tau)
\mathop{{\sum}}\limits_{l_r=0}^{\frac{n}{(j-1)^\gamma}-r} \sqrt{l_r+1/2}\: \varphi_{l_r}(t) \,
\sum_{k=\frac{n}{(j-1)^\gamma}}^{n}\, \sqrt{k+1/2} \:\langle f, \varphi_{k,j}\:
 \rangle B^r_{k} ,
$$
$$
\triangle_{22}(t,\tau)= 2^r \sum_{j=2}^{(n/r)^{1/\gamma}+1} \varphi_{j}(\tau)
\mathop{{\sum}}\limits_{l_r=\frac{n}{(j-1)^\gamma}-r+1}^{n-r} \sqrt{l_r+1/2}\: \varphi_{l_r}(t) \,
\sum_{k=l_r+r}^{n}\, \sqrt{k+1/2} \:\langle f, \varphi_{k,j}\:
 \rangle B^r_{k} ,
$$
$$
\triangle_{23}(t,\tau)= 2^r \sum_{j=(n/r)^{1/\gamma}+2}^{\infty} \varphi_{j}(\tau)
\mathop{{\sum}}\limits_{l_r=0}^{n-r} \sqrt{l_r+1/2}\: \varphi_{l_r}(t) \, \sum_{k=l_r+r}^{n}\, \sqrt{k+1/2} \:\langle
f, \varphi_{k,j}\:
 \rangle B^r_{k} .
$$
Let's first estimate the norm of $\triangle_{21}$. Again using the formula (\ref{Muller}) and H\"{o}lder inequality, we
get:
$$
\|\triangle_{21}\|_{L_2}^2\leq c \sum_{j=2}^{(n/r)^{1/\gamma}+1} j^{-2\mu_2}
\mathop{{\sum}}\limits_{l_r=0}^{\frac{n}{(j-1)^\gamma}-r} (l_r+1/2)
$$
$$
\times \left( \sum_{k=\frac{n}{(j-1)^\gamma}}^{n} k^{s\mu_1} j^{s\mu_2} |\langle f, \varphi_{kj}\rangle
|^s\right)^{2/s} \left( \sum_{k=\frac{n}{(j-1)^\gamma}}^{n} k^{-(\mu_1-2r+3/2)s/(s-1)} \right)^{2(s-1)/s}
$$
$$
\leq c  \, n^{-2(\mu_1-2r+1/2)+\frac{2(s-1)}{s}} \sum_{j=2}^{(n/r)^{1/\gamma}+1}
j^{-2\mu_2+2\gamma(\mu_1-2r+1/2)-\frac{2\gamma(s-1)}{s}} \left( \sum_{k=\frac{n}{(j-1)^\gamma}}^{n} k^{s\mu_1}
j^{s\mu_2} |\langle f, \varphi_{k,j}\rangle |^s\right)^{2/s}
$$
$$
\leq c \|f\|_{s,\overline{\mu}}^2 \, n^{-2(\mu_1-2r+1/s-1/2)}.
$$
Then
$$
\|\triangle_{21}\|_{L_2}\leq c \|f\|_{s,\overline{\mu}} \, n^{-\mu_1+2r-1/s+1/2}.
$$
Further we can bound the norm of $\triangle_{22}(t,\tau):$
$$
\|\triangle_{22}\|_{L_2}^2 \leq c \, \sum_{j=2}^{(n/r)^{1/\gamma}+1} j^{-2\mu_2} \left( \sum_{k=r}^{n} k^{s\mu_1}
j^{s\mu_2} |\langle f, \varphi_{k,j}\rangle |^s\right)^{2/s} \,
\mathop{{\sum}}\limits_{l_r=\frac{n}{(j-1)^\gamma}-r+1}^{n-r} (l_r+1/2)^{-2(\mu_1-2r+1)+\frac{2(s-1)}{s}}
$$
$$
\leq c  \, n^{-2(\mu_1-2r+1/s-1/2)} \sum_{j=2}^{(n/r)^{1/\gamma}+1} j^{-2\mu_2+2\gamma(\mu_1-2r+1/s-1/2)} \left(
\sum_{k=r}^{n} k^{s\mu_1} j^{s\mu_2} |\langle f, \varphi_{kj}\rangle |^s\right)^{2/s}
$$
$$
= c \|f\|_{s,\overline{\mu}}^2 \, n^{-2(\mu_1-2r+1/s-1/2)}.
$$
Similar to the previous we find
$$
\|\triangle_{23}\|_{L_2}^2 \leq c \, \sum_{j=(n/r)^{1/\gamma}+2}^{\infty} j^{-2\mu_2} \left( \sum_{k=r}^{n} k^{s\mu_1}
j^{s\mu_2} |\langle f, \varphi_{k,j}\rangle |^s\right)^{2/s} \, \mathop{{\sum}}\limits_{l_r=0}^{n-r}
(l_r+1/2)^{-2(\mu_1-2r+1)+\frac{2(s-1)}{s}}
$$
$$
\leq c \, \sum_{j=(n/r)^{1/\gamma}+2}^{\infty} j^{-2\mu_2} \left( \sum_{k=r}^{n} k^{s\mu_1} j^{s\mu_2} |\langle f,
\varphi_{k,j}\rangle |^s\right)^{2/s}
$$
$$
\leq c  \|f\|_{s,\overline{\mu}}^2 \, n^{-(2\mu_2-\frac{s-2}{s})/\gamma} \leq c \|f\|_{s,\overline{\mu}}^2 \,
n^{-2(\mu_1-2r+1/s-1/2)}.
$$

Combining the estimates for $\triangle_{21}(t,\tau)$, $\triangle_{22}(t,\tau)$ and $\triangle_{23}(t,\tau)$ we obtain
$$
\|\triangle_{2}\|_{L_2}
\leq c\, \|f\|_{s,\overline{\mu}} \, n^{-\mu_1+2r-1/s+1/2}.
$$

The combination of (\ref{Bound_err}) and bounds for the norms of $\triangle_{1}$, $\triangle_{2}$ makes it possible to
establish the desired inequality.

Finding the estimate of the norm $(\ref{Bound_err})$ in the case of $s=2$ is similar to the previous one. \vspace{0.1in}

\vskip -4mm

  ${}$ \ \ \ \ \ \ \ \ \ \ \ \ \ \ \ \ \ \ \ \ \ \ \ \ \ \ \ \ \ \ \ \ \ \ \ \ \ \
  \ \ \ \ \ \ \ \ \ \ \ \ \ \ \ \ \ \ \ \ \ \
  \ \ \ \ \ \ \ \ \ \ \ \ \ \ \ \ \ \ \ \ \ \
  \ \ \ \ \ \ \ \ \ \ \ \ \ \ \ \ \ \ \ \ \ \
  \ \ \ \ \ \ \ \ \ \ \ \ \ \ \ \ \ \ \ \ \ \
  \ \ \ \ \ \ \ \ \ \ \ \ \ \ \ \ \ \ \ \ \ \
  \ \ $\Box$

In a similar way, the following statement is established.

\begin{lemma}\label{lemma_Bound1+1}
Let $f$ belongs to the class $L^{\overline{\mu}}_{s,2}$, $1\leq s< 2$, $\mu_1>2r-1/s+1/2$, $\mu_2\geq
\mu_1-2r+1/s-1/2$.
 Then for $1\leq \gamma \leq \frac{\mu_2}{\mu_1-2r+1/s-1/2}$ it holds
$$
  \|f^{(r,0)}-\mathcal{D}_{n,\gamma}^{(r,0)} f\|_{L_2}\leq c\|f\|_{s,\overline{\mu}}^{s/2} n^{-\mu_1+2r-1/s+1/2} .
$$
\end{lemma}

 \vspace{0.1in}

The following statement contains an estimate for the second difference from the right-hand side of (\ref{fullError}) in
the metric of $L_2$.

\begin{lemma}\label{lemma_BoundPertHC}
Let the condition (\ref{perturbation}) be satisfied for $1\leq p \leq \infty$.  Then for any function $f\in
L_2(Q)$ and $\gamma$ it holds 
 $$
 \|\mathcal{D}^{(r,0)}_{n,\gamma} f - \mathcal{D}^{(r,0)}_{n,\gamma} f^\delta\|_{L_2} \leq c\delta n^{2r-1/p+1/2} .
 $$
\end{lemma}
\textit{Proof.} Let's write down the representation
$$\mathcal{D}^{(r,0)}_{n,\gamma} f(t,\tau) - \mathcal{D}^{(r,0)}_{n,\gamma} f^\delta(t,\tau)=
\Big(\sum_{k=r}^{n} \sum_{j=0}^{(n/k)^{1/\gamma}} \langle f-f^\delta, \varphi_{k,j}\rangle
\varphi'_k(t)\varphi_j(\tau)\Big)^{(r-1)}_t.
$$
Using the formula (\ref{Muller}) we get
$$\mathcal{D}^{(r,0)}_{n,\gamma} f(t,\tau) - \mathcal{D}^{(r,0)}_{n,\gamma} f^\delta(t,\tau)= 2^r \sum_{k=r}^{n}\,
\sum_{j=0}^{(n/k)^{1/\gamma}} \sqrt{k+1/2}\, \langle f-f^\delta, \varphi_{k,j}\rangle \varphi_j(\tau)
$$
$$
\times \mathop{{\sum}}\limits_{l_1=r-1}^{k-1}  (l_1+1/2)\: \mathop{{\sum}}\limits_{l_2=r-2}^{l_1-1} (l_2+1/2)\:
\ldots \mathop{{\sum}}\limits_{l_{r-1}=1}^{l_{r-2}-1} (l_{r-1}+1/2)\:
 \mathop{{\sum}}\limits_{l_{r}=0}^{l_{r-1}-1} \sqrt{l_{r}+1/2}\: \varphi_{l_r}(t) .
 $$
Further, we change the order of summation and get
$$
\mathcal{D}^{(r,0)}_{n,\gamma} f(t,\tau) - \mathcal{D}^{(r,0)}_{n,\gamma} f^\delta(t,\tau) = 2^r\,
\sum_{j=0}^{(n/r)^{1/\gamma}} \, \varphi_j(\tau) \mathop{{\sum}}\limits_{l_r=0}^{n/\underline{j}^\gamma-r}
\sqrt{l_r+1/2}\, \varphi_{l_r}(t)
$$
$$
\times \sum_{k=l_r+r}^{n/\underline{j}^\gamma} \sqrt{k+1/2}\, \langle f-f^\delta, \varphi_{k,j}\: \rangle B^r_{k} .
$$
Let $1<p<\infty$ first. Then, using the H\"{o}lder inequality and the estimate (\ref{axular1}), we find
$$
\|\mathcal{D}^{(r,0)}_{n,\gamma} f(t,\tau) - \mathcal{D}^{(r,0)}_{n,\gamma} f^\delta\|_{L_2}^2 \leq c\,
\sum_{j=0}^{(n/r)^{1/\gamma}} \, \mathop{{\sum}}\limits_{l_r=0}^{n/\underline{j}^\gamma-r} (l_r+1/2)
\left(\sum_{k=l_r+r}^{n/\underline{j}^\gamma} |\langle f-f^\delta, \varphi_{k,j}\: \rangle | k^{2r-3/2}\right)^2
$$
$$
\leq c \delta^2\, \sum_{j=0}^{(n/r)^{1/\gamma}} \, \mathop{{\sum}}\limits_{l_r=0}^{n/\underline{j}^\gamma-r} (l_r+1/2)
\left(\sum_{k=l_r+r}^{n/\underline{j}^\gamma} k^{(2r-3/2)p/(p-1)}\right)^{2(p-1)/p}
$$
$$
\leq c\delta^2 n^{2(2r-3/2)+2(p-1)/p} \sum_{j=0}^{(n/r)^{1/\gamma}} \, \underline{j}^{-2\gamma(2r-3/2)-2\gamma(p-1)/p}
\mathop{{\sum}}\limits_{l_r=0}^{n/\underline{j}^\gamma-r} (l_r+1/2)
$$
$$
 = c \delta^2 n^{2(2r-1/2)+2(p-1)/p} ,
$$
which was required to prove.

In the case of $p=1$ and $p=\infty$, the assertion of Lemma is proved similarly. \vspace{0.1in}

\vskip -4mm

  ${}$ \ \ \ \ \ \ \ \ \ \ \ \ \ \ \ \ \ \ \ \ \ \ \ \ \ \ \ \ \ \ \ \ \ \ \ \ \ \
  \ \ \ \ \ \ \ \ \ \ \ \ \ \ \ \ \ \ \ \ \ \
  \ \ \ \ \ \ \ \ \ \ \ \ \ \ \ \ \ \ \ \ \ \
  \ \ \ \ \ \ \ \ \ \ \ \ \ \ \ \ \ \ \ \ \ \
  \ \ \ \ \ \ \ \ \ \ \ \ \ \ \ \ \ \ \ \ \ \
  \ \ \ \ \ \ \ \ \ \ \ \ \ \ \ \ \ \ \ \ \ \
  \ \ $\Box$

The combination of Lemmas \ref{lemma_BoundErrHC} and \ref{lemma_BoundPertHC} gives
\begin{theorem} \label{Th1}
Let $f\in L^{\overline{\mu}}_{s,2}$, $2\leq s< \infty$, $\mu_1>2r-1/s+1/2$, $\mu_2>\mu_1-2r$, and let the condition
(\ref{perturbation}) be satisfied for $1\leq p \leq \infty$.
 Then for $n\asymp \delta^{-\frac{1}{\mu_1-1/p+1/s}}$ and $1\leq \gamma < \frac{\mu_2+1/s-1/2}{\mu_1-2r+1/s-1/2}$ it holds
  $$
   \|f^{(r,0)} - \mathcal{D}^{(r,0)}_{n,\gamma} f^\delta\|_{L_2} \leq c \delta^{\frac{\mu_1-2r+1/s-1/2}{\mu_1-1/p+1/s}} .
  $$
\end{theorem}

\vskip 2mm

\begin{corollary} \label{Cor1}
 \rm In the considered problem, the truncation method  $\mathcal{D}^{(r,0)}_{n,\gamma}$ (\ref{ModVer})
 achieves the accuracy
$O\left(\delta^{\frac{\mu_1-2r+1/s-1/2}{\mu_1-1/p+1/s}}\right)$
 on the class $L^{\overline{\mu}}_{s,2}$, $2\leq s< \infty$, $\mu_1>2r-1/s+1/2$, $\mu_2>\mu_1-2r$, and requires
 $$
 \mathrm{card}(\Gamma_{n,\gamma}) \asymp \left\{
\begin{array}{cl}
n \asymp \delta^{-\frac{1}{\mu_1-1/p+1/s}}, & \mbox{ if }\ 1 < \gamma < \frac{\mu_2+1/s-1/2}{\mu_1-2r+1/s-1/2}, \\
n\ \ln n \asymp \delta^{-\frac{1}{\mu_1-1/p+1/s}} \ln \frac{1}{\delta}, & \mbox{ if }\ \gamma=1,
\end{array}
\right.
 $$
perturbed Fourier-Legendre coefficients.
\end{corollary}

\vskip 2mm

The combination of Lemmas \ref{lemma_Bound1+1} and \ref{lemma_BoundPertHC} gives
\begin{theorem} \label{Th1+1}
Let $f\in L^{\overline{\mu}}_{s,2}$, $1\leq s< 2$, $\mu_1>2r-1/s+1/2$, $\mu_2\geq \mu_1-2r+1/s-1/2$, and let the
condition (\ref{perturbation}) be satisfied for $1\leq p \leq \infty$.
 Then for $n\asymp \delta^{-\frac{1}{\mu_1-1/p+1/s}}$ and $1\leq \gamma \leq \frac{\mu_2}{\mu_1-2r+1/s-1/2}$ it holds
  $$
   \|f^{(r,0)} - \mathcal{D}^{(r,0)}_{n,\gamma} f^\delta\|_{L_2} \leq c \delta^{\frac{\mu_1-2r+1/s-1/2}{\mu_1-1/p+1/s}} .
  $$
\end{theorem}

\vskip 2mm

\begin{corollary} \label{Cor1+1}
 \rm In the considered problem, the truncation method  $\mathcal{D}^{(r,0)}_{n,\gamma}$ (\ref{ModVer})
 achieves the accuracy
$O\left(\delta^{\frac{\mu_1-2r+1/s-1/2}{\mu_1-1/p+1/s}}\right)$
 on the class $L^{\overline{\mu}}_{s,2}$, $1\leq s< 2$, $\mu_1>2r-1/s+1/2$, $\mu_2\geq \mu_1-2r+1/s-1/2$, and requires
 $$
 \mathrm{card}(\Gamma_{n,\gamma}) \asymp \left\{
\begin{array}{cl}
n \asymp \delta^{-\frac{1}{\mu_1-1/p+1/s}}, & \mbox{ if }\ 1 < \gamma \leq \frac{\mu_2}{\mu_1-2r+1/s-1/2}, \\
n\ \ln n \asymp \delta^{-\frac{1}{\mu_1-1/p+1/s}} \ln \frac{1}{\delta}, & \mbox{ if }\ \gamma=1,
\end{array}
\right.
 $$
perturbed Fourier-Legendre coefficients.
\end{corollary}

\vskip 2mm

\section{Truncation Method. Error estimate in the metric of $C$}

Now we have to bound the error of (\ref{ModVer}) in the metric of $C$. An upper estimate for the norm of the difference
(\ref{Bound_err}) is contained in the following statement.

\begin{lemma}\label{lemma_BoundErrHCC}
Let $f\in L^{\overline{\mu}}_{s,2}$, $1\leq s< \infty$, $\mu_1>2r-1/s+3/2$, $\mu_2>\mu_1-2r$.
 Then for $1\leq \gamma < \frac{\mu_2+1/s-3/2}{\mu_1-2r+1/s-3/2}$ it holds
$$
  \|f^{(r,0)}-\mathcal{D}^{(r,0)}_{n,\gamma} f\|_{C}\leq c\|f\|_{s,\overline{\mu}} \, n^{-\mu_1+2r-1/s+3/2} .
$$
\end{lemma}

\textit{Proof.}
Using (\ref{triangle_{11}}) and (\ref{axular1}),  we get for $1< s< \infty$
$$
\|\triangle_{11}\|_{C} \leq c \sum_{j=0}^{\infty} \underline{j}^{-\mu_2+1/2} \mathop{{\sum}}\limits_{l_r=0}^{n-r+1}
(l_r+1/2)
$$
$$
\times \left( \sum_{k=n+1}^{\infty} k^{s\mu_1} \underline{j}^{s\mu_2} |\langle f, \varphi_{k,j}\rangle |^s\right)^{1/s}
\left( \sum_{k=n+1}^{\infty} k^{-(\mu_1-2r+3/2)s/(s-1)} \right)^{(s-1)/s}
$$
$$
\leq c  \|f\|_{s,\overline{\mu}} \, n^{-\mu_1+2r+1/2+\frac{s-1}{s}} \Big(\sum_{j=0}^{\infty}
\underline{j}^{-\frac{s}{s-1}(\mu_2-1/2)}\Big)^{(s-1)/s} = c \|f\|_{s,\overline{\mu}} \, n^{-\mu_1+2r-1/s+3/2} .
$$
Moreover, from (\ref{triangle_{12}}) it follows
$$
\|\triangle_{12}\|_{C} \leq c n^{-\mu_1+2r+1/2+\frac{s-1}{s}} \, \sum_{j=0}^{\infty} \underline{j}^{-\mu_2+1/2} \left(
\sum_{k=r}^{\infty} k^{s\mu_1} \underline{j}^{s\mu_2} |\langle f, \varphi_{k,j}\rangle |^s\right)^{1/s}
 $$
$$
\leq c \|f\|_{s,\overline{\mu}} \, n^{-\mu_1+2r-1/s+3/2} .
$$
Thus we get
$$
\|\triangle_{1}\|_{C} \leq c \|f\|_{s,\overline{\mu}} \, n^{-\mu_1+2r-1/s+3/2} .
$$
Similarly, we find
$$
\|\triangle_{21}\|_{C} \leq c \sum_{j=2}^{(n/r)^{1/\gamma}+1} j^{-\mu_2+1/2}
\mathop{{\sum}}\limits_{l_r=0}^{\frac{n}{(j-1)^\gamma}-r} (l_r+1/2)
$$
$$
\times \left( \sum_{k=\frac{n}{(j-1)^\gamma}}^{n} k^{s\mu_1} j^{s\mu_2} |\langle f, \varphi_{k,j}\rangle
|^s\right)^{1/s} \left( \sum_{k=\frac{n}{(j-1)^\gamma}}^{n} k^{-(\mu_1-2r+3/2)s/(s-1)} \right)^{(s-1)/s}
$$
$$
\leq c  \, n^{-\mu_1+2r+1/2+\frac{s-1}{s}} \sum_{j=2}^{(n/r)^{1/\gamma}+1} j^{-\mu_2+1/2+\gamma(\mu_1-2r+1/s-3/2)}
\left( \sum_{k=\frac{n}{(j-1)^\gamma}}^{n} k^{s\mu_1} j^{s\mu_2} |\langle f, \varphi_{kj}\rangle |^s\right)^{1/s}
$$
$$
\leq c \|f\|_{s,\overline{\mu}} \, n^{-\mu_1+2r-1/s+3/2} ,
$$
$$
\|\triangle_{22}\|_{C} \leq c \, n^{-\mu_1+2r-1/s+3/2} \, \sum_{j=2}^{(n/r)^{1/\gamma}+1}
j^{-\mu_2+1/2+\gamma(\mu_1-2r+1/s-3/2)} \left( \sum_{k=r}^{n} k^{s\mu_1} j^{s\mu_2} |\langle f, \varphi_{k,j}\rangle
|^s\right)^{1/s}
$$
$$
\leq c \|f\|_{s,\overline{\mu}} \, n^{-\mu_1+2r-1/s+3/2} ,
$$
$$
\|\triangle_{23}\|_{C} \leq c \, \sum_{j=(n/r)^{1/\gamma}+2}^{\infty} j^{-\mu_2+1/2} \left( \sum_{k=r}^{n} k^{s\mu_1}
j^{s\mu_2} |\langle f, \varphi_{k,j}\rangle |^s\right)^{1/s} \mathop{{\sum}}\limits_{l_r=0}^{n-r}
(l_r+1/2)^{-\mu_1+2r-1/2+\frac{s-1}{s}}
$$
$$
\leq c \|f\|_{s,\overline{\mu}} \, n^{-(\mu_2+1/s-3/2)/\gamma}
\leq c \|f\|_{s,\overline{\mu}} \, n^{-\mu_1+2r-1/s+3/2}
.
$$
Combining the estimates obtained above, we have
$$
\|\triangle_{2}\|_{C} \leq c \|f\|_{s,\overline{\mu}} \, n^{-\mu_1+2r-1/s+3/2} .
$$

Substituting estimates for the norms of $\triangle_{1}$, $\triangle_{2}$ into the relation (\ref{Bound_err}) allows to
establish the desired inequality.

Finding the estimate of the norm $(\ref{Bound_err})$ in the case of $s=1$ is similar to the previous one. \vskip -4mm

  ${}$ \ \ \ \ \ \ \ \ \ \ \ \ \ \ \ \ \ \ \ \ \ \ \ \ \ \ \ \ \ \ \ \ \ \ \ \ \ \
  \ \ \ \ \ \ \ \ \ \ \ \ \ \ \ \ \ \ \ \ \ \
  \ \ \ \ \ \ \ \ \ \ \ \ \ \ \ \ \ \ \ \ \ \
  \ \ \ \ \ \ \ \ \ \ \ \ \ \ \ \ \ \ \ \ \ \
  \ \ \ \ \ \ \ \ \ \ \ \ \ \ \ \ \ \ \ \ \ \
  \ \ \ \ \ \ \ \ \ \ \ \ \ \ \ \ \ \ \ \ \ \
  \ \ $\Box$

The following statement contains an estimate for the second difference from the right-hand side of (\ref{fullError}) in
the metric of $C$.
\begin{lemma}\label{lemma_BoundPertHCC}
Let the condition (\ref{perturbation}) be satisfied for $1\leq p \leq \infty$.  Then for any function $f\in
L_2(Q)$ and $\gamma \geq 1$  it holds
 $$
 \|\mathcal{D}^{(r,0)}_{n,\gamma} f - \mathcal{D}^{(r,0)}_{n,\gamma} f^\delta\|_{C} \leq c \delta n^{2r-1/p+3/2} .
 $$
\end{lemma}
\textit{Proof.} Let $1< p<\infty$ first. Then, using the H\"{o}lder inequality and the estimate (\ref{axular1}), we find
$$
\|\mathcal{D}^{(r,0)}_{n,\gamma} f - \mathcal{D}^{(r,0)}_{n,\gamma} f^\delta\|_{C} \leq c\,
\sum_{j=0}^{(n/r)^{1/\gamma}} \, \mathop{{\sum}}\limits_{l_r=0}^{n/\underline{j}^\gamma-r} (l_r+1/2)
\sum_{k=l_r+r}^{n/\underline{j}^\gamma} |\langle f-f^\delta, \varphi_{kj}\: \rangle | k^{2r-3/2}
$$
$$
\leq c \delta \, \sum_{j=0}^{(n/r)^{1/\gamma}} \, \mathop{{\sum}}\limits_{l_r=0}^{n/\underline{j}^\gamma-r} (l_r+1/2)
\left(\sum_{k=l_r+r}^{n/\underline{j}^\gamma} k^{(2r-3/2)p/(p-1)}\right)^{(p-1)/p}
$$
$$
\leq c\delta n^{2r-3/2+(p-1)/p} \sum_{j=0}^{(n/r)^{1/\gamma}} \, \underline{j}^{-\gamma(2r-3/2)-\gamma(p-1)/p}
\mathop{{\sum}}\limits_{l_r=0}^{n/\underline{j}^\gamma-r} (l_r+1/2)
$$
$$
 \leq c \delta n^{2r-1/p+3/2} ,
$$
which was required to prove.

In the case of $p=1$ and $p=\infty$, the assertion of Lemma is proved similarly. \vspace{0.1in}

\vskip -4mm

  ${}$ \ \ \ \ \ \ \ \ \ \ \ \ \ \ \ \ \ \ \ \ \ \ \ \ \ \ \ \ \ \ \ \ \ \ \ \ \ \
  \ \ \ \ \ \ \ \ \ \ \ \ \ \ \ \ \ \ \ \ \ \
  \ \ \ \ \ \ \ \ \ \ \ \ \ \ \ \ \ \ \ \ \ \
  \ \ \ \ \ \ \ \ \ \ \ \ \ \ \ \ \ \ \ \ \ \
  \ \ \ \ \ \ \ \ \ \ \ \ \ \ \ \ \ \ \ \ \ \
  \ \ \ \ \ \ \ \ \ \ \ \ \ \ \ \ \ \ \ \ \ \
  \ \ $\Box$

The combination of Lemmas \ref{lemma_BoundErrHCC} and \ref{lemma_BoundPertHCC} gives
\begin{theorem} \label{Th2}
Let $f\in L^{\overline{\mu}}_{s,2}$, $1\leq s< \infty$, $\mu_1>2r-1/s+3/2$, $\mu_2>\mu_1-2r$, and let the condition
(\ref{perturbation}) be satisfied for $1\leq p \leq \infty$.
 Then for $n\asymp \delta^{-\frac{1}{\mu_1-1/p+1/s}}$ and $1\leq \gamma < \frac{\mu_2+1/s-3/2}{\mu_1-2r+1/s-3/2}$ it holds
  $$
   \|f^{(r,0)} - \mathcal{D}^{(r,0)}_{n,\gamma} f^\delta\|_{C} \leq c \delta^{\frac{\mu_1-2r+1/s-3/2}{\mu_1-1/p+1/s}} .
  $$
\end{theorem}

\vskip 2mm

\begin{corollary} \label{Cor2}
 \rm In the considered problem, the truncation method  $\mathcal{D}^{(r,0)}_{n,\gamma}$ (\ref{ModVer})
 achieves the accuracy
$O\left(\delta^{\frac{\mu_1-2r+1/s-3/2}{\mu_1-1/p+1/s}}\right)$
 on the class $L^{\overline{\mu}}_{s,2}$, $1\leq s< \infty$, $\mu_1>2r-1/s+3/2$, $\mu_2>\mu_1-2r$, and requires
 $$
 \mathrm{card}(\Gamma_{n,\gamma}) \asymp \left\{
\begin{array}{cl}
n \asymp \delta^{-\frac{1}{\mu_1-1/p+1/s}}, & \mbox{ if }\ 1 < \gamma < \frac{\mu_2+1/s-3/2}{\mu_1-2r+1/s-3/2}, \\
n\ \ln n \asymp \delta^{-\frac{1}{\mu_1-1/p+1/s}} \ln \frac{1}{\delta}, & \mbox{ if }\ \gamma=1,
\end{array}
\right.
 $$
perturbed Fourier-Legendre coefficients.
\end{corollary}

\vskip 2mm

\vskip 2mm

\begin{remark} \label{Rem3}
\rm Previously (see \cite{Sol_Stas_UMZ2022}) the problem of the minimal Galerkin information radius for the problem of numerical differentiation of functions from  $L^{\overline{\mu}}_{s,2}$ was studied in the situation where $r=1$ and $p=s=2$.
Thus, the results of Theorems \ref{Th1}, \ref{Th1+1} and \ref{Th2}  generalize the investigations
of \cite{Sol_Stas_UMZ2022} to the case of arbitrary $r,p,s$.

\end{remark}

\vskip 2mm
\section{Minimal radius of Galerkin information}

Let us turn to find sharp estimates (in the power scale) for the minimal radius. First, we establish a lower
estimate for the quantity $R_{N,\delta}^{(r,0)}(L^{\overline{\mu}}_{s,2}, C, \ell_p)$. We fix an arbitrarily chosen
domain
 $\hat{\Omega}$, $\card(\hat{\Omega})\leq N$, of the coordinate plane $[r,\infty)\times[0,\infty)$ and build an auxiliary function
 $$
 f_1(t,\tau)  = \widetilde{c} \, \bigg( \varphi_0(t) \varphi_0(\tau) \,
 + \,  N^{-\mu_1-1/s} \varphi_1(\tau) \mathop{{\sum}\, '}\limits_{k=N+r}^{3N+r}
 \varphi_k(t) \bigg) ,
 $$
 where the sum $\mathop{{\sum}\, '}\limits_{k=N+r}^{3N+r}$ is taken over any $N$ pairwise distinct functions $\varphi_k(t)$ such that
 $N+r\leq k\leq 3N+r$ and $(k,1)\notin \hat{\Omega}$. Obviously, there is at least one set of such functions.


Now we estimate the norm of $f_1$ in the space metric of $L^{\overline{\mu}}_{s,2}$:
 $$
 \|f_1\|_{s,\overline{\mu}} = {\widetilde{c}} \, \bigg( 1
 +  N^{-s\mu_1-1}  \,  \mathop{{\sum}\, '}\limits_{k=N+r}^{3N+r}
 k^{s\mu_1} \bigg)^{1/s} \leq
{\widetilde{c}}    \,  \bigg( 1
 + 4^{s\mu_1} \bigg)^{1/s} .
 $$
 Whence it follows that to satisfy the condition $\|f_1\|_{s,\overline{\mu}}\leq 1$ it suffices to take
 $$
 {\widetilde{c}} =  \bigg( 1
 + 4^{s\mu_1} \bigg)^{-1/s} .
 $$

Next, we take another function from the class $L^{\overline{\mu}}_{s,2}$:
 $$
 f_2(t,\tau) = \widetilde{c} \, \varphi_0(t) \varphi_0(\tau)  .
 $$

Let us find a lower bound for the quantity $\|f_1^{(r,0)}-f_2^{(r,0)}\|_{C}$. For this we need formulas
$$
\varphi_1(1) = \sqrt{3/2} , \quad f_2^{(r,0)}(t,\tau) \equiv 0 ,
$$
 $$
 f_1^{(r,0)}(t,\tau) = \widetilde{c} \, {N^{-\mu_1-1/s}} \,
\varphi_1(\tau) \mathop{{\sum}\, '}\limits_{k=N+r}^{3N+r} \varphi_k^{(r)}(t)
$$
$$
= 2^r\, \widetilde{c} \, {N^{-\mu_1-1/s}} \, \varphi_1(\tau) \mathop{{\sum}\, '}\limits_{k=N+r}^{3N+r} \sqrt{k+1/2}
\mathop{{\sum}}\limits_{l_1=r-1}^{k-1} (l_1+1/2) \mathop{{\sum}}\limits_{l_2=r-2}^{l_1-1} (l_2+1/2)
$$
\begin{equation}  \label{f_1^(r,r)}
\ldots \mathop{{\sum}}\limits_{l_{r-1}=1}^{l_{r-2}-1} (l_{r-1}+1/2) \mathop{{\sum}}\limits_{l_{r}=0}^{l_{r-1}-1}
\sqrt{l_r+1/2}\, \varphi_{l_r}(t) .
\end{equation}
We note that in the right-hand side of (\ref{f_1^(r,0)}) only terms with odd indexes $l_1+k, l_2+l_1,..., l_r+l_{r-1}$ take part.

It is easy to see that
$$
 \|f_1^{(r,0)}-f_2^{(r,0)}\|_{C}
 \geq  |f_1^{(r,0)}(1,1)|
 \geq  \overline{c} \, N^{-\mu_1+2r-1/s+3/2}  ,
 $$
where
 $$
\overline{c} = \frac{\sqrt{3/2}}{2^{r} r!}\, \widetilde{c} .
 $$
Since for any $1\leq p\leq \infty$ it holds true
 $$
 \|\overline{f}_1-\overline{f}_2\|_{\ell_p}
 = \widetilde{c} \, N^{-\mu_1-1/s+1/p} ,
 $$
 then in the case of $N^{-\mu_1-1/s+1/p}\leq \delta/{\widetilde{c}}$  under $\delta$-perturbations of the functions $f_1$ and $f_2$
 can be considered
$$
f^\delta_1(t,\tau) = f_2(t,\tau), \qquad f^\delta_2(t,\tau) = f_1(t,\tau) .
$$

Let us find  the upper bound  for $\|f_1^{(r,0)}-f_2^{(r,0)}\|_{C}$.
 Taking into account the relation $G(\hat{\Omega},\overline{f}_1^{\delta})=G(\hat{\Omega},\overline{f}_2^{\delta})$,
for any $\psi^{(r,0)}(\hat{\Omega})\in\Psi(\hat{\Omega})$ we find
 $$
 \|f_1^{(r,0)}-f_2^{(r,0)}\|_{C}
 \leq \|f_1^{(r,0)}-\psi^{(r,0)}(G(\hat{\Omega},\overline{f}_1^{\delta}))\|_{C}
 + \|f_2^{(r,0)}-\psi^{(r,0)}(G(\hat{\Omega},\overline{f}_2^{\delta}))\|_{C}
 \leq
 $$
 $$
 \leq 2 \, \sup\limits_{f\in L^{\overline{\mu}}_{s,2}, \, \|f\|_{s,\overline{\mu}}\leq 1}
 \ \sup\limits_{\overline{f^\delta}: \, (\ref{perturbation})}
 \| f^{(r,0)} - \psi^{(r,0)}(G(\hat{\Omega},\overline{f}^{\delta})) \|_C
 =: 2 \, \varepsilon_{\delta}(L^{\overline{\mu}}_{s,2}, \psi^{(r,0)}(\hat{\Omega}), C, \ell_p) .
 $$

That is
 $$
 \varepsilon_{\delta}(L^{\overline{\mu}}_{s,2}, \psi^{(r,0)}(\hat{\Omega}), C, \ell_p)
 \geq \frac{\overline{c}}{2} \, N^{-\mu_1+2r-1/s+3/2} .
 $$

From the fact that the domain $\hat{\Omega}$
 and the algorithm $\psi^{(r,0)}(\hat{\Omega})\in\Psi(\hat{\Omega})$  are arbitrary, it follows
 $$
R_{N,\delta}^{(r,0)}(L^{\overline{\mu}}_{s,2}, C, \ell_p)
 \geq \frac{\overline{c}}{2} \, N^{-\mu_1+2r-1/s+3/2} .
 $$

\vskip 2mm

Thus, the following statement is proved.

\begin{theorem} \label{Th5.1}
Let $1\leq s< \infty$, $\mu_1>2r-1/s+3/2$, $1\leq p \leq \infty$.
 Then for any $N\geq\Big(\delta/\widetilde{c}\Big)^{-1/(\mu_1+1/s-1/p)}$ it holds
 $$
R_{N,\delta}^{(r,0)}(L^{\overline{\mu}}_{s,2}, C, \ell_p)
 \geq \frac{\overline{c}}{2} \, N^{-\mu_1+2r-1/s+3/2} .
 $$
\end{theorem}

The following assertion contains sharp estimates (in the power scale) for the minimal radius in the uniform metric.

\begin{theorem} \label{Th5.2}
Let $1\leq s< \infty$, $\mu_1>2r-1/s+3/2$, $\mu_2>\mu_1-2r$, $1\leq p \leq \infty$.

a)  If $N \asymp \delta^{-\frac{1}{\mu_1-1/p+1/s}}$, then it holds
$$
R_{N,\delta}^{(r,0)}(L^{\overline{\mu}}_{s,2}, C, \ell_p)
 \asymp N^{-\mu_1+2r-1/s+3/2} \asymp \delta^{\frac{\mu_1-2r+1/s-3/2}{\mu_1-1/p+1/s}} .
$$
The indicated order-optimal estimates are implemented by the method $\mathcal{D}^{(r,0)}_{n,\gamma}$ (\ref{ModVer}) for
$n\asymp \delta^{-\frac{1}{\mu_1-1/p+1/s}}$ and $1<\gamma<\frac{\mu_2+1/s-3/2}{\mu_1-2r+1/s-3/2}$.

b)  If $N \asymp \delta^{-\frac{1}{\mu_1-1/p+1/s}} \ln \frac{1}{\delta}$, then it holds
$$
N^{-\mu_1+2r-1/s+3/2} \preceq R_{N,\delta}^{(r,0)}(L^{\overline{\mu}}_{s,2}, C, \ell_p) \preceq (N/\ln
N)^{-\mu_1+2r-1/s+3/2}
$$
or
$$
\delta^{\frac{\mu_1-2r+1/s-3/2}{\mu_1-1/p+1/s}} \ln^{-\mu_1+2r-1/s+3/2} \frac{1}{\delta} \preceq
R_{N,\delta}^{(r,0)}(L^{\overline{\mu}}_{s,2}, C, \ell_p) \preceq \delta^{\frac{\mu_1-2r+1/s-3/2}{\mu_1-1/p+1/s}} .
$$
The upper bounds are implemented by the method $\mathcal{D}^{(r,0)}_{n,\gamma}$ (\ref{ModVer})  for  $n\asymp
\delta^{-\frac{1}{\mu_1-1/p+1/s}}$ and $\gamma=1$.
\end{theorem}

 \bf Proof.
 \rm
The upper bounds for  $R_{N,\delta}^{(r,0)}(L^{\overline{\mu}}_{s,2}, C, \ell_p)$ follow from Theorem \ref{Th2}.
The lower bound is found in Theorem  \ref{Th5.1}.

\vskip -4mm

  ${}$ \ \ \ \ \ \ \ \ \ \ \ \ \ \ \ \ \ \ \ \ \ \ \ \ \ \ \ \ \ \ \ \ \ \ \ \ \ \
  \ \ \ \ \ \ \ \ \ \ \ \ \ \ \ \ \ \ \ \ \ \
  \ \ \ \ \ \ \ \ \ \ \ \ \ \ \ \ \ \ \ \ \ \
  \ \ \ \ \ \ \ \ \ \ \ \ \ \ \ \ \ \ \ \ \ \
  \ \ \ \ \ \ \ \ \ \ \ \ \ \ \ \ \ \ \ \ \ \
  \ \ \ \ \ \ \ \ \ \ \ \ \ \ \ \ \ \ \ \ \ \
  \ \ $\Box$

 \rm

\vskip 2mm

Let's turn to estimate the minimal radius in the integral metric.

\begin{theorem} \label{Th5.4}
Let $1\leq s< \infty$, $\mu_1>2r-1/s+1/2$, $1\leq p \leq \infty$. Then for any
$N\geq\Big(\delta/\widetilde{c}\Big)^{-1/(\mu_1+1/s-1/p)}$ it holds
$$
R_{N,\delta}^{(r,0)}(L^{\overline{\mu}}_{s,2}, L_2, \ell_p) \geq \overline{\overline{c}} \, N^{-\mu_1+2r-1/s+1/2} ,
$$
 where $\overline{\overline{c}} = \frac{\widetilde{c}}{2^{3r-1}(r-1)!}$ .
\end{theorem}

 \rm

 \bf Proof \rm The proof of Theorem \ref{Th5.4} almost completely coincides with the proof of Theorem \ref{Th5.1},
 including the form of the auxiliary functions
 $f_1$, $f_1^{\delta}$, $f_2$, $f_2^{\delta}$.
 The only difference is in the lower estimate of the norm of the difference
 $f_1^{(r,0)}-f_2^{(r,0)}$. Changing the order of summation in  (\ref{f_1^(r,r)}) yields to the representation
$$
 f_1^{(r,0)}(t,\tau) = 2^r\, \widetilde{c} \, {N^{-\mu_1-1/s}} \, \varphi_1(\tau)
 \Big(\mathop{{\sum}}\limits_{l_r=0}^{N} \sqrt{l_r+1/2}\, \varphi_{l_r}(t) \mathop{{\sum}\,
'}\limits_{k=N+r}^{3N+r} \sqrt{k+1/2}
$$
$$
+ \mathop{{\sum}}\limits_{l_r=N+1}^{3N} \sqrt{l_r+1/2}\, \varphi_{l_r}(t) \mathop{{\sum}\, '}\limits_{k=l_r+r}^{3N+r}
\sqrt{k+1/2}\Big) \mathop{{\sum}}\limits_{l_1=l_r+r-1}^{k-1} (l_1+1/2) \mathop{{\sum}}\limits_{l_2=l_r+r-2}^{l_1-1}
(l_2+1/2) \ldots \mathop{{\sum}}\limits_{l_{r-1}=l_r+1}^{l_{r-2}-1} (l_{r-1}+1/2)  .
$$
Then it is easy to see
 $$
 \|f_1^{(r,0)}-f_2^{(r,0)}\|_{L_2}^2 = \|f_1^{(r,0)}\|_{L_2}^2
 $$
$$
\geq 4^r\, \widetilde{c}^2 \, {N^{-2\mu_1-2/s}}
 \mathop{{\sum}}\limits_{l_r=0}^{N/2} (l_r+1/2)\, \left( \mathop{{\sum}\,
'}\limits_{k=N+r}^{3N+r} \sqrt{k+1/2} \frac{(k-N/2-r+1)^{2r-2}}{4^{r-1}(r-1)!}\right)^2
$$
$$
  \geq \frac{\widetilde{c}^2}{4^{3r-4}((r-1)!)^2} \, {N^{-2\mu_1-2/s}} \cdot N^{4r-1}
 \mathop{{\sum}}\limits_{l_r=0}^{N/2} (l_r+1/2) .
 $$
That is, it holds
 $$
 \|f_1^{(r,0)}-f_2^{(r,0)}\|_{L_2}
 \geq \frac{\widetilde{c}}{2^{3r-2}(r-1)!} \,  N^{-\mu_1+2r-1/s+1/2} .
 $$
 Whence we obtain the relation
 $$
 \varepsilon_{\delta}(L^{\overline{\mu}}_{s,2}, \psi^{(r,0)}(\hat{\Omega}), L_2, \ell_p)
 \geq \overline{\overline{c}} \, N^{-\mu_1+2r-1/s+1/2}
 $$
is true for any $N\geq\Big(\delta/\widetilde{c}\Big)^{-1/(\mu_1+1/s-1/p)}$. From the fact that the domain
$\hat{\Omega}$ and the algorithm $\psi^{(r,0)}(\hat{\Omega})\in\Psi(\hat{\Omega})$  are arbitrary, it follows that
 $$
R_{N,\delta}^{(r,0)}(L^{\overline{\mu}}_{s,2}, L_2, \ell_p)
 \geq \overline{\overline{c}} \, N^{-\mu_1+2r-1/s+1/2} .
 $$
Thus, Theorem \ref{Th5.4} is proved.

\vskip -4mm

  ${}$ \ \ \ \ \ \ \ \ \ \ \ \ \ \ \ \ \ \ \ \ \ \ \ \ \ \ \ \ \ \ \ \ \ \ \ \ \ \
  \ \ \ \ \ \ \ \ \ \ \ \ \ \ \ \ \ \ \ \ \ \
  \ \ \ \ \ \ \ \ \ \ \ \ \ \ \ \ \ \ \ \ \ \
  \ \ \ \ \ \ \ \ \ \ \ \ \ \ \ \ \ \ \ \ \ \
  \ \ \ \ \ \ \ \ \ \ \ \ \ \ \ \ \ \ \ \ \ \
  \ \ \ \ \ \ \ \ \ \ \ \ \ \ \ \ \ \ \ \ \ \
  \ \ $\Box$

\vskip 2mm

The next two statements contain order estimates of the minimal radius in the integral metric.
 \rm

\begin{theorem} \label{Th5.5}
Let $2\leq s< \infty$, $\mu_1>2r-1/s+1/2$, $\mu_2>\mu_1-2r$, $1\leq p \leq \infty$.

a)  If $N \asymp \delta^{-\frac{1}{\mu_1-1/p+1/s}}$, then it holds true
$$
R_{N,\delta}^{(r,0)}(L^{\overline{\mu}}_{s,2}, L_2, \ell_p)
 \asymp N^{-\mu_1+2r-1/s+1/2} \asymp \delta^{\frac{\mu_1-2r+1/s-1/2}{\mu_1-1/p+1/s}} .
$$
The indicated order-optimal estimates are implemented by the method $\mathcal{D}^{(r,0)}_{n,\gamma}$ (\ref{ModVer}) for
$n\asymp \delta^{-\frac{1}{\mu_1-1/p+1/s}}$ and $1<\gamma<\frac{\mu_2+1/s-1/2}{\mu_1-2r+1/s-1/2}$.

b)  If $N \asymp \delta^{-\frac{1}{\mu_1-1/p+1/s}} \ln \frac{1}{\delta}$, then it holds true
$$
N^{-\mu_1+2r-1/s+1/2} \preceq R_{N,\delta}^{(r,0)}(L^{\overline{\mu}}_{s,2}, L_2, \ell_p) \preceq (N/\ln
N)^{-\mu_1+2r-1/s+1/2}
$$
or
$$
\delta^{\frac{\mu_1-2r+1/s-1/2}{\mu_1-1/p+1/s}} \ln^{-\mu_1+2r-1/s+1/2} \frac{1}{\delta} \preceq
R_{N,\delta}^{(r,0)}(L^{\overline{\mu}}_{s,2}, L_2, \ell_p) \preceq \delta^{\frac{\mu_1-2r+1/s-1/2}{\mu_1-1/p+1/s}} .
$$
The upper bounds are implemented by the method $\mathcal{D}^{(r,0)}_{n,\gamma}$ (\ref{ModVer})  for  $n\asymp
\delta^{-\frac{1}{\mu_1-1/p+1/s}}$ and $\gamma=1$.
\end{theorem}

{\textbf{\textit{Proof.}}
 \rm
The upper estimates for $R_{N,\delta}^{(r,0)}(L^{\overline{\mu}}_{s,2}, L_2, \ell_p)$ follow from Theorem \ref{Th1 }.
The lower estimate of the minimal radius is found in Theorem \ref{Th5.4}.

\vskip -4mm

  ${}$ \ \ \ \ \ \ \ \ \ \ \ \ \ \ \ \ \ \ \ \ \ \ \ \ \ \ \ \ \ \ \ \ \ \ \ \ \ \
  \ \ \ \ \ \ \ \ \ \ \ \ \ \ \ \ \ \ \ \ \ \
  \ \ \ \ \ \ \ \ \ \ \ \ \ \ \ \ \ \ \ \ \ \
  \ \ \ \ \ \ \ \ \ \ \ \ \ \ \ \ \ \ \ \ \ \
  \ \ \ \ \ \ \ \ \ \ \ \ \ \ \ \ \ \ \ \ \ \
  \ \ \ \ \ \ \ \ \ \ \ \ \ \ \ \ \ \ \ \ \ \
  \ \ $\Box$

\vskip 2mm

The combination of Theorems \ref{Th1+1} and \ref{Th5.4} gives the next assertion.

\begin{theorem} \label{Th5.5+1}
Let $1\leq s< 2$, $\mu_1>2r-1/s+1/2$, $\mu_2\geq \mu_1-2r+1/s-1/2$, $1\leq p \leq \infty$.

a)  If $N \asymp \delta^{-\frac{1}{\mu_1-1/p+1/s}}$, then it holds
$$
R_{N,\delta}^{(r,0)}(L^{\overline{\mu}}_{s,2}, L_2, \ell_p)
 \asymp N^{-\mu_1+2r-1/s+1/2} \asymp \delta^{\frac{\mu_1-2r+1/s-1/2}{\mu_1-1/p+1/s}} .
$$
The indicated order-optimal estimates are implemented by the method $\mathcal{D}^{(r,0)}_{n,\gamma}$ (\ref{ModVer}) for
$n\asymp \delta^{-\frac{1}{\mu_1-1/p+1/s}}$ and $1<\gamma\leq \frac{\mu_2}{\mu_1-2r+1/s-1/2}$.

b)  If $N \asymp \delta^{-\frac{1}{\mu_1-1/p+1/s}} \ln \frac{1}{\delta}$, then it holds
$$
N^{-\mu_1+2r-1/s+1/2} \preceq R_{N,\delta}^{(r,0)}(L^{\overline{\mu}}_{s,2}, L_2, \ell_p) \preceq (N/\ln
N)^{-\mu_1+2r-1/s+1/2}
$$
or
$$
\delta^{\frac{\mu_1-2r+1/s-1/2}{\mu_1-1/p+1/s}} \ln^{-\mu_1+2r-1/s+1/2} \frac{1}{\delta} \preceq
R_{N,\delta}^{(r,0)}(L^{\overline{\mu}}_{s,2}, L_2, \ell_p) \preceq \delta^{\frac{\mu_1-2r+1/s-1/2}{\mu_1-1/p+1/s}} .
$$
The upper bounds are implemented by the method $\mathcal{D}^{(r,0)}_{n,\gamma}$ (\ref{ModVer})  for  $n\asymp
\delta^{-\frac{1}{\mu_1-1/p+1/s}}$ and $\gamma=1$.
\end{theorem}

\vskip 2mm

Finally, we consider the problem of optimal recovering the derivative $f^{(0,r)}$ in the sense of quantity
 $$
 R^{(0,r)}_{N,\delta} (L^{\overline{\mu}}_{s,2}, X, \ell_2) = \inf\limits_{\Omega: \, \card(\Omega)\leq N}
 \ \, \inf\limits_{\psi^{(0,r)}\in\Psi(\Omega)} \varepsilon_{\delta}(L^{\overline{\mu}}_{s,2}, \psi^{(0,r)}(\Omega), X, \ell_p) .
 $$
where
 $$
 \varepsilon_{\delta}(L^{\overline{\mu}}_{s,2}, \psi^{(0,r)}(\Omega), X, \ell_p)
 = \sup\limits_{f\in L^{\overline{\mu}}_{s,2}, \, \|f\|_{s,\overline{\mu}}\leq 1}
 \ \, \sup\limits_{\overline{f^\delta}: \, (\ref{perturbation})}
 \| f^{(0,r)} - \psi^{(0,r)}(G(\Omega, \overline{f^\delta})) \|_X
 $$
is the error of the algorithm $\psi^{(0,r)}$ on the class $L^{\overline{\mu}}_{s,2}$, $X=L_2$ or $X=C$,
and by the $r$-th partial derivative $f^{(0,r)}$ of the function $f\in L^{\overline{\mu}}_{s,2}$ we mean the series
 $$ 
 f^{(0,r)}(t,\tau) = \sum\limits_{k=0}^{\infty} \sum\limits_{j=r}^{\infty}
 \langle f , \varphi_{k,j} \rangle \, \varphi_{k}(t) \, \varphi_{j}^{(r)}(\tau) .
 $$ 
.

Let
 \begin{equation}\label{D_n_f_delta_(0,r)}
\mathcal{D}^{(0,r)}_{n,\gamma}f^{\delta}(t,\tau) = \sum\limits_{j\geq r,\, k^{\gamma}j\leq n}
 \langle f_{\delta} , \varphi_{k,j} \rangle \, \varphi_{k}(t) \, \varphi_{j}^{(r)}(\tau) ,
 \quad \gamma \geq 1 ,
 \end{equation}
 be the proposed variant of the truncation method to recovering $f^{(0,r)}$, where
the hyperbolic cross is taken as the domain $\Omega$ of the coordinate plane $[0,\infty)\times [r,\infty)$
 $$
\Gamma_{n,\gamma} = \{(k,j):\, k^{\gamma}j\leq n,\quad 0 \leq k \leq n^{1/\gamma}, \quad r\leq j \leq n\} .
 $$

Analogously to Theorems \ref{Th5.2}, \ref{Th5.5}, \ref{Th5.5+1}, the following statements can be obtained.

\vskip 2mm

\begin{theorem} \label{Th5.2.1}
Let $1\leq s<\infty$, $\mu_2>2r-1/s+3/2$, $\mu_1>\mu_2-2r$, $1\leq p \leq \infty$.

a) If $N \asymp \delta^{-\frac{1}{\mu_2-1/p+1/s}}$, then it holds
$$
R_{N,\delta}^{(0,r)}(L^{\overline{\mu}}_{s,2}, C, \ell_p)
 \asymp N^{-\mu_2+2r-1/s+3/2} \asymp \delta^{\frac{\mu_2-2r+1/s-3/2}{\mu_2-1/p+1/s}} .
$$
The indicated order-optimal estimates are implemented by the method $\mathcal{D}^{(0,r)}_{n,\gamma}$
(\ref{D_n_f_delta_(0,r)}) for $n\asymp \delta^{-\frac{1}{\mu_2-1/p+1/s}}$ and
$1<\gamma<\frac{\mu_1+1/s-3/2}{\mu_2-2r+1/s-3/2}$.

b)  If $N \asymp \delta^{-\frac{1}{\mu_2-1/p+1/s}} \ln \frac{1}{\delta}$, then it holds
$$
N^{-\mu_2+2r-1/s+3/2} \preceq R_{N,\delta}^{(0,r)}(L^{\overline{\mu}}_{s,2}, C, \ell_p) \preceq (N/\ln
N)^{-\mu_2+2r-1/s+3/2}
$$
or
$$
\delta^{\frac{\mu_2-2r+1/s-3/2}{\mu_2-1/p+1/s}} \ln^{-\mu_2+2r-1/s+3/2} \frac{1}{\delta} \preceq
R_{N,\delta}^{(0,r)}(L^{\overline{\mu}}_{s,2}, C, \ell_p) \preceq \delta^{\frac{\mu_2-2r+1/s-3/2}{\mu_2-1/p+1/s}} .
$$
The upper bounds are implemented by the method $\mathcal{D}^{(0,r)}_{n,\gamma}$ (\ref{D_n_f_delta_(0,r)})  for $n\asymp
\delta^{-\frac{1}{\mu_2-1/p+1/s}}$ and $\gamma=1$.
\end{theorem}

\vskip 2mm

\begin{theorem} \label{Th5.5.1}
Let $2\leq s< \infty$, $\mu_2>2r-1/s+1/2$, $\mu_1>\mu_2-2r$, $1\leq p \leq \infty$.

a)  If $N \asymp \delta^{-\frac{1}{\mu_2-1/p+1/s}}$, then it holds
$$
R_{N,\delta}^{(0,r)}(L^{\overline{\mu}}_{s,2}, L_2, \ell_p)
 \asymp N^{-\mu_2+2r-1/s+1/2} \asymp \delta^{\frac{\mu_2-2r+1/s-1/2}{\mu_2-1/p+1/s}} .
$$
The indicated order-optimal estimates are implemented by the method $\mathcal{D}^{(0,r)}_{n,\gamma}$
(\ref{D_n_f_delta_(0,r)}) for $n\asymp \delta^{-\frac{1}{\mu_2-1/p+1/s}}$ and
$1<\gamma<\frac{\mu_1+1/s-1/2}{\mu_2-2r+1/s-1/2}$.

b)  If $N \asymp \delta^{-\frac{1}{\mu_2-1/p+1/s}} \ln \frac{1}{\delta}$, then it holds
$$
N^{-\mu_2+2r-1/s+1/2} \preceq R_{N,\delta}^{(0,r)}(L^{\overline{\mu}}_{s,2}, L_2, \ell_p) \preceq (N/\ln
N)^{-\mu_2+2r-1/s+1/2}
$$
or
$$
\delta^{\frac{\mu_2-2r+1/s-1/2}{\mu_2-1/p+1/s}} \ln^{-\mu_2+2r-1/s+1/2} \frac{1}{\delta} \preceq
R_{N,\delta}^{(0,r)}(L^{\overline{\mu}}_{s,2}, L_2, \ell_p) \preceq \delta^{\frac{\mu_2-2r+1/s-1/2}{\mu_2-1/p+1/s}} .
$$
The upper bounds are implemented by the method $\mathcal{D}^{(0,r)}_{n,\gamma}$ (\ref{D_n_f_delta_(0,r)})  for $n\asymp
\delta^{-\frac{1}{\mu_2-1/p+1/s}}$ and $\gamma=1$.
\end{theorem}

\begin{theorem} \label{Th5.5.1+1}
Let $1\leq s< 2$, $\mu_2>2r-1/s+1/2$, $\mu_1\geq \mu_2-2r+1/s-1/2$, $1\leq p \leq \infty$.

a)  If $N \asymp \delta^{-\frac{1}{\mu_2-1/p+1/s}}$, then it holds
$$
R_{N,\delta}^{(0,r)}(L^{\overline{\mu}}_{s,2}, L_2, \ell_p)
 \asymp N^{-\mu_2+2r-1/s+1/2} \asymp \delta^{\frac{\mu_2-2r+1/s-1/2}{\mu_2-1/p+1/s}} .
$$
The indicated order-optimal estimates are implemented by the method $\mathcal{D}^{(0,r)}_{n,\gamma}$
(\ref{D_n_f_delta_(0,r)}) for $n\asymp \delta^{-\frac{1}{\mu_2-1/p+1/s}}$ and $1<\gamma\leq
\frac{\mu_1}{\mu_2-2r+1/s-1/2}$.

b)  If $N \asymp \delta^{-\frac{1}{\mu_2-1/p+1/s}} \ln \frac{1}{\delta}$, then it holds
$$
N^{-\mu_2+2r-1/s+1/2} \preceq R_{N,\delta}^{(0,r)}(L^{\overline{\mu}}_{s,2}, L_2, \ell_p) \preceq (N/\ln
N)^{-\mu_2+2r-1/s+1/2}
$$
or
$$
\delta^{\frac{\mu_2-2r+1/s-1/2}{\mu_2-1/p+1/s}} \ln^{-\mu_2+2r-1/s+1/2} \frac{1}{\delta} \preceq
R_{N,\delta}^{(0,r)}(L^{\overline{\mu}}_{s,2}, L_2, \ell_p) \preceq \delta^{\frac{\mu_2-2r+1/s-1/2}{\mu_2-1/p+1/s}} .
$$
The upper bounds are implemented by the method $\mathcal{D}^{(0,r)}_{n,\gamma}$ (\ref{D_n_f_delta_(0,r)})  for $n\asymp
\delta^{-\frac{1}{\mu_2-1/p+1/s}}$ and $\gamma=1$.
\end{theorem}

\section{Computational experiments}

To demonstrate the effectiveness of the proposed method $\mathcal{D}^{(r,0)}_{n,\gamma}$ (\ref{ModVer}) some numerical experiments were carried out. The calculations were performed on a computer with a 4-core Intel Core i5 processor and 16 GB memory in the mathematical modeling environment MATLAB 2022a.


\subsection{Example 1.}

We consider the function  $F(t,\tau) = f(t)f(\tau)/C$,  where $C=947$  and
$$
f(t)=\left\{
\begin{array}{cl}
	-1/8 t^2+1/12 t^4-1/25 t^5+1/38 t^7-1/108 t^8,   & -1\leq t<0 ,
	\\\\
	-1/8 t^2+1/12 t^4-1/25 t^5+1/102 t^7-1/198 t^8, &  0\leq t \leq 1.
\end{array}
\right.
$$	
It is easy to see that for $\mu=5,6$ we have $\|F\|_{2,\mu} \approx 1$ and $\|F^{(2,0)}\|_{L_2}\approx 10^{-5}.$

The simulation of the noise in the input data was done in two different ways:

\begin{itemize}
	\item a random noise adds to the values of the Fourier-Legendre coefficients. The noise is generated by the $\mbox{randn(size}({\cal F})) \delta$ command, where randn and size are standard functions of the MATLAB system, and ${\cal F}$ is a matrix for exact values of the Fourier-Legendre coefficients;
	
	\item the values of the Fourier-Legendre coefficients recover by the quadrature trapezoid formula on a uniform grid with a step $h$   so that condition (\ref{perturbation}) is satisfied for a given $\delta$.

	%
\end{itemize}

Numerical experiments were carried out for the following error levels: $\delta= 10^{-7}, 10^{-8}, 10^{-9}$. Tables \ref{tbl1} and \ref{tbl2} show the results of numerical calculations for the approximation  of $F^{(2,0)}$ by the truncation method (\ref{ModVer})
with two different types of noise (random data noise and trapezoid formula errors).  ErrorL2 and ErrorC columns contain the recovery accuracy in the $L_2$ and $C$ metrics respectively, the $n$ and $ \card(\Gamma_{n})$ columns indicate the highest degree of Legendre polynomial and the number of Fourier-Legendre coefficients involved, respectively.
Also in Table \ref{tbl2}, $h$ means the step size in the quadrature formula.

Graphs \ref{Fig2} and \ref{Fig1} show the exact derivative $F^{(2,0)}$ and its approximations constructed on data with random noise and with noise,
generated by the trapezoid formula, respectively.

%

\begin{table}[h!]
	\centering
	\caption{ The results of recovering derivative $F^{(2,0)}$ for random noise }
	\label{tbl1}
	\begin{tabular}{|c|c|c|c|}
		\hline
		$\delta$ & $10^{-7}$ & $10^{-8}$ & $10^{-9}$  \\ \hline
		ErrorL2        &  $2,1 \cdot 10^{-6} $       &  $8,6 \cdot 10^{-7} $  &  $3,4  \cdot 10^{-7} $  \\ \hline
		ErrorC        &  $1,8 \cdot 10^{-5} $       &  $5,58 \cdot 10^{-6} $  &  $1,37\cdot10^{-6} $  \\ \hline
		$n$  &  16      & 25    & 28    \\ \hline
	\end{tabular}
\end{table}
\begin{table}[ht]
	\centering
	\caption{ The results of recovering derivative $F^{(2,0)}$ for noise from quadrature formula }
	\label{tbl2}
	\begin{tabular}{|c|c|c|c|}
		\hline
		$\delta$ & $10^{-7}$ & $10^{-8}$ & $10^{-9}$  \\ \hline
		ErrorL2        &  $1,9 \cdot 10^{-6} $       &  $1,6 \cdot 10^{-6} $  &  $4,5  \cdot 10^{-7} $  \\ \hline
		ErrorC        &  $1,8 \cdot 10^{-5} $       &  $1,38 \cdot 10^{-5} $  &  $3,14 \cdot10^{-6} $  \\ \hline
		n  &  16      & 22    & 28    \\ \hline
		h  &  $1\cdot 10^{-4}$    & $8\cdot10^{-5} $  & $ 4\cdot10^{-5}$    \\ \hline
	\end{tabular}
\end{table}

\begin{figure}[h!]
	\begin{minipage}[h]{0.5\linewidth}
		\center{\includegraphics[width=1\linewidth]{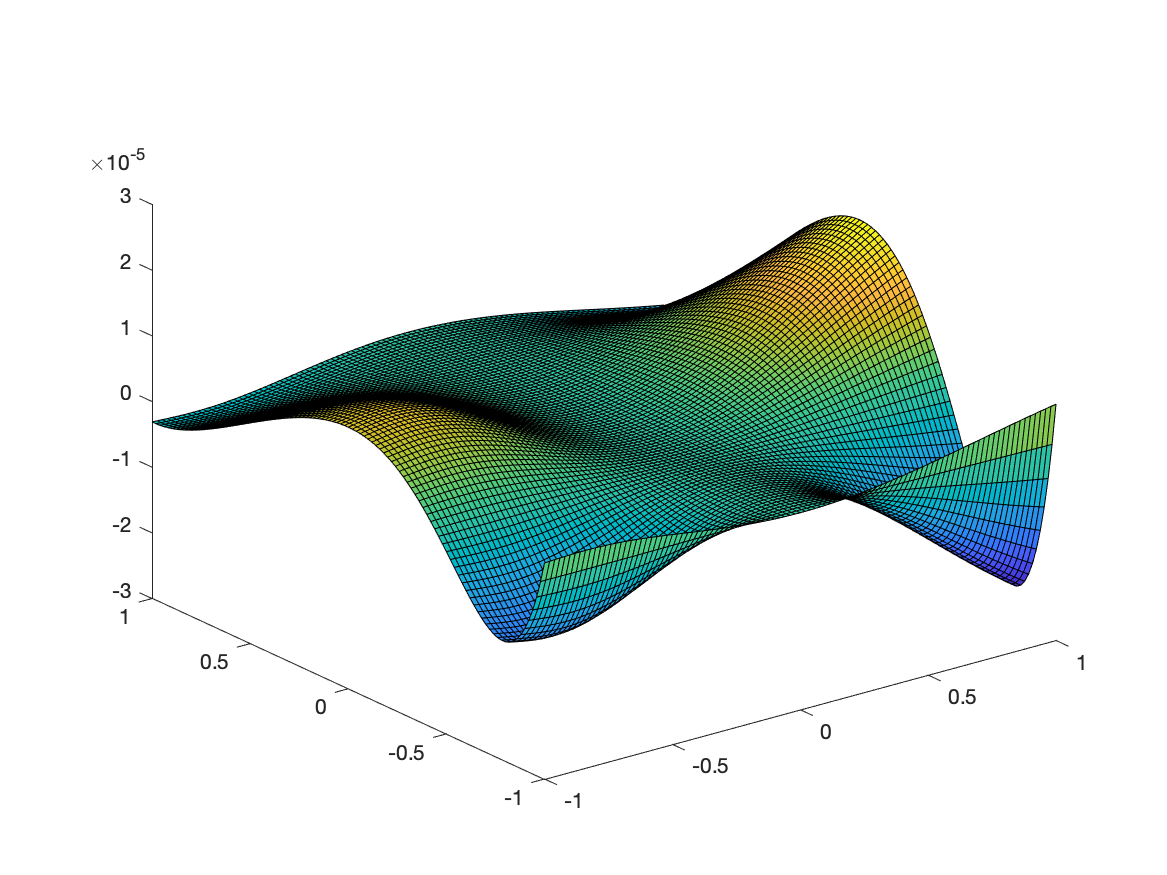} \\ a)}
	\end{minipage}
	\begin{minipage}[h]{0.5\linewidth}
		\center{\includegraphics[width=1\linewidth]{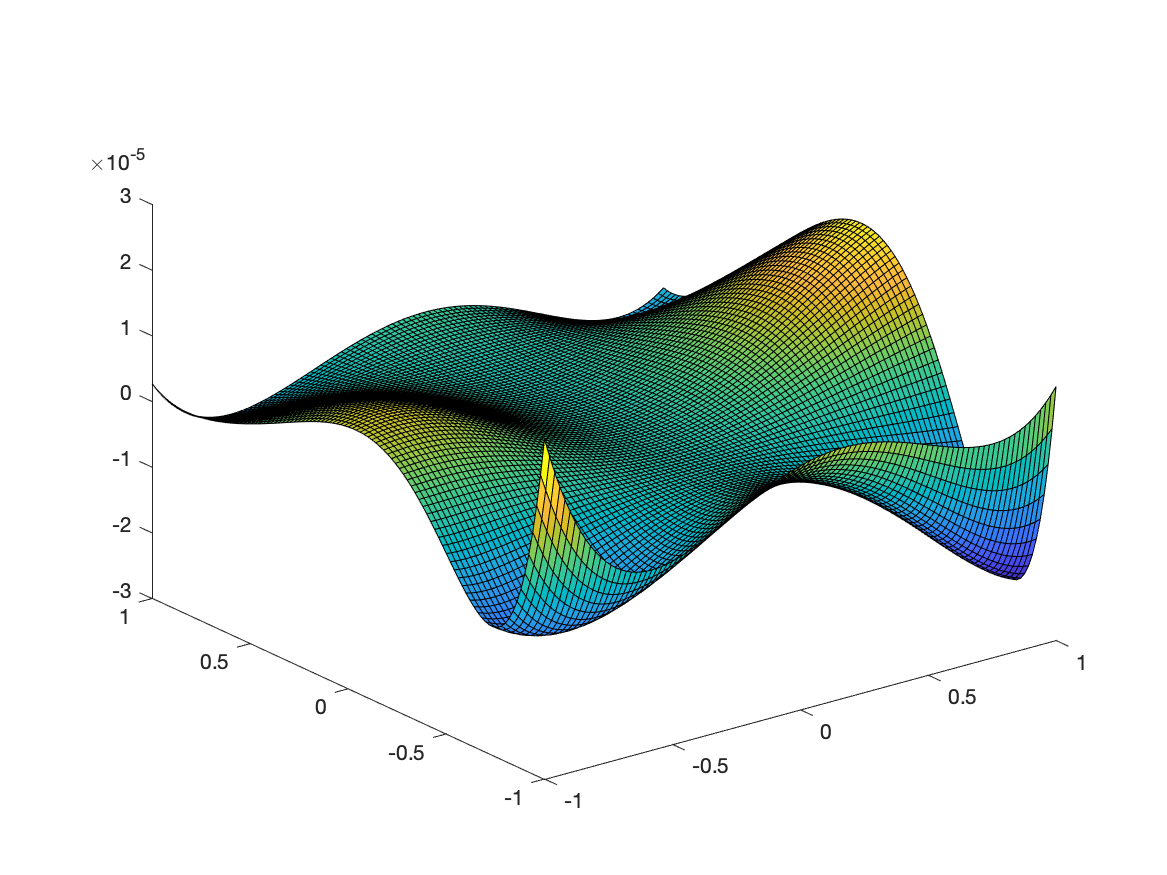} \\ b)}
	\end{minipage}
	\hfill
	\begin{minipage}[h]{0.5\linewidth}
		\center{\includegraphics[width=1\linewidth]{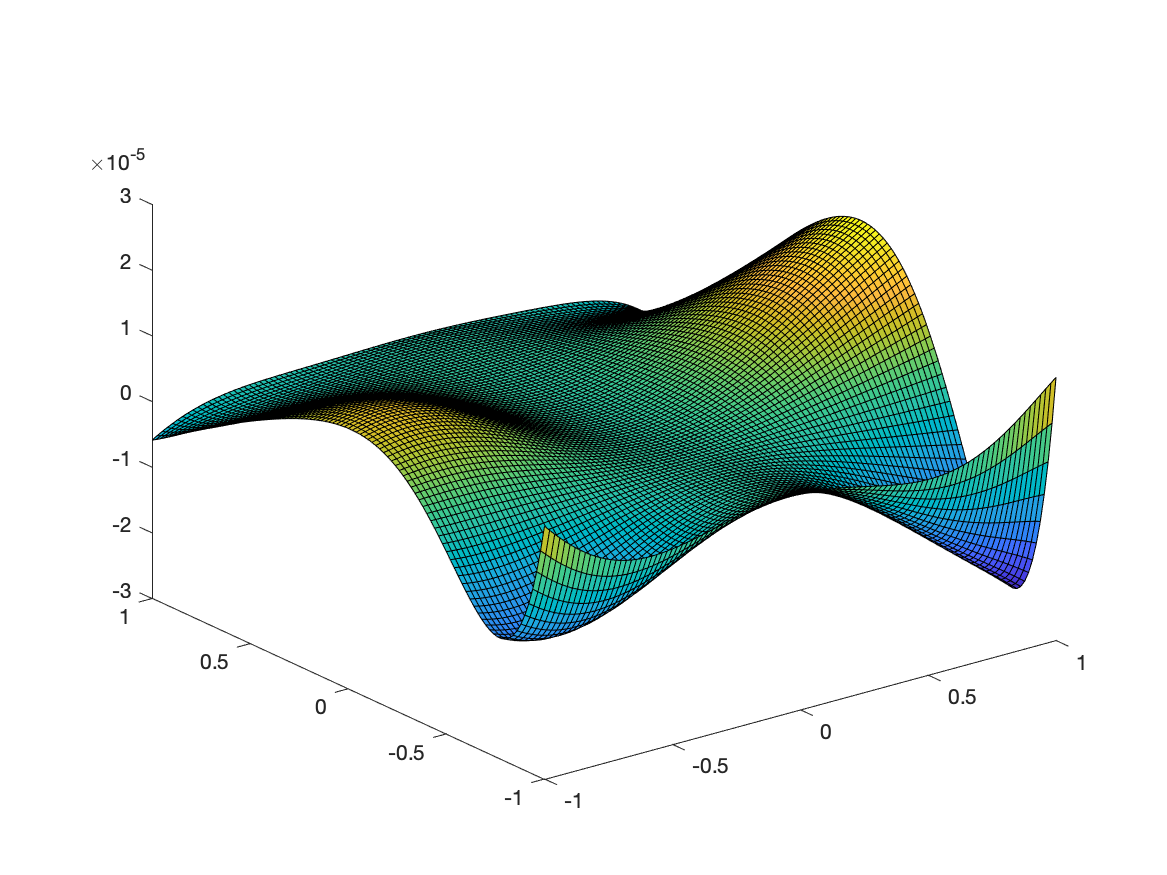} \\ c)}
	\end{minipage}
	\begin{minipage}[h]{0.5\linewidth}
		\center{\includegraphics[width=1\linewidth]{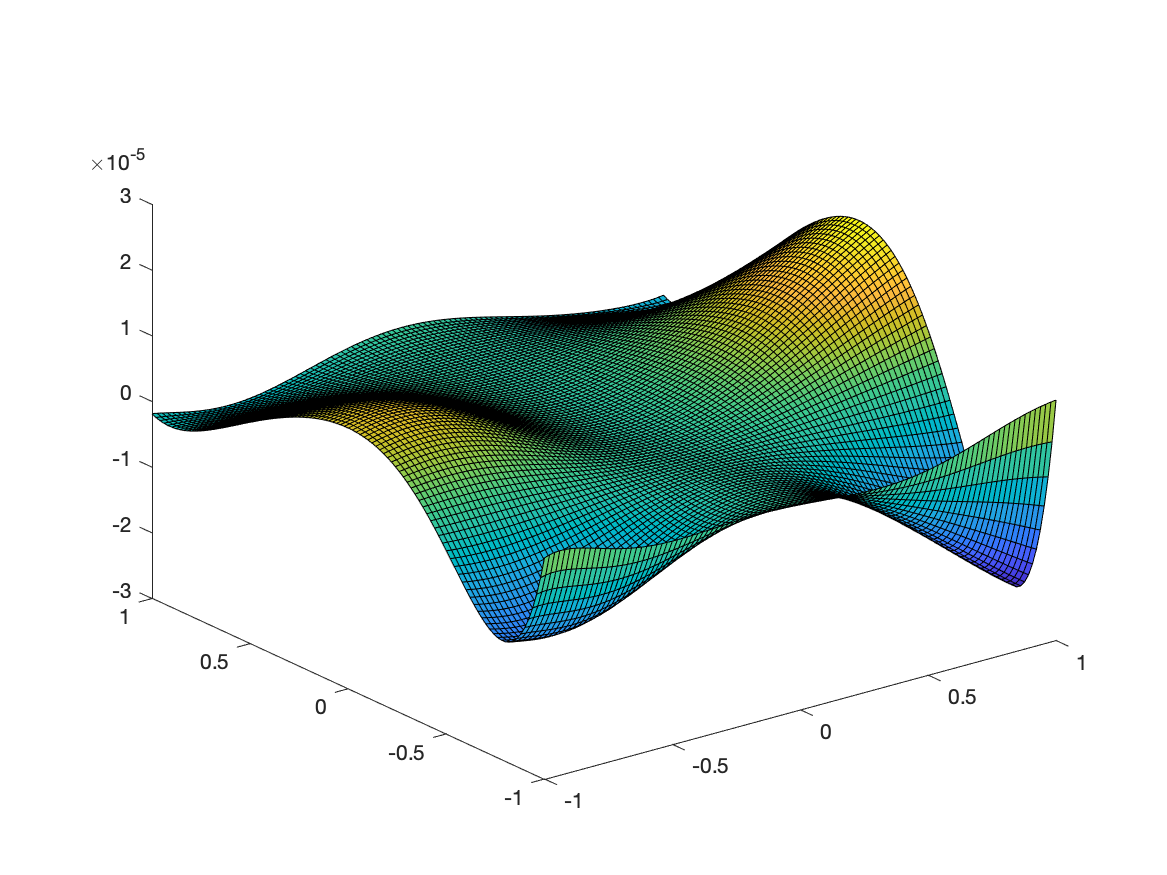} \\ d)}
	\end{minipage}
	\caption{Recovery of the derivative $F^{(2,0)}$  with  random noise in the  input data . The exact derivative $F^{(2,0)}$ (Fig. а) ); approximation to $F^{(2,0)}$ for    $\delta= 10^{-7}$ (Fig.  b) );    for $\delta= 10^{-8}$ (Fig. с) ) and $\delta= 10^{-9}$ (Fig.  d) ), }
	\label{Fig2}
\end{figure}

As can be seen from the graphs and tables above, for both types of noise, the truncation  method (\ref{ModVer})
gives the same order of accuracy for recovering $F^{(2,0)}$ .
At the same time, applying the quadrature formula expands the area of using the proposed method in computational problems, especially in the situation when the input data are given in the form of a set of function values at the grid nodes.

\subsection{Example 2}

Let us test the method (\ref{ModVer}) on  the function $ F(t,\tau)= f_1(t)f_2(\tau)/C, $ where $f_2(t)=2\cos(\pi t)$ and
$$
f_1(t)=\left\{
\begin{array}{cl}
	-1/8 t^2+1/12 t^4-1/25 t^5+1/38 t^7-1/108 t^8,   & -1\leq t<0 ,
	\\\\
	-1/8 t^2+1/12 t^4-1/25 t^5+1/102 t^7-1/198 t^8, &  0\leq t \leq 1.
\end{array}
\right.
$$	
Let's put $\mu=5.4$. It is easy to check that $\|F\|_{2,\mu} \approx 1$ and $\|F^{(2,0)}\|_{L_2} \approx 10^{-5}$, if $C=26318$.
The Fourier-Legendre coefficients are calculated using the quadrature trapezoid formula for $h = 8\cdot 10^{-5}, 2\cdot10^{-5}, 8\cdot10^{-6}$, which in turn according to formula (\ref{perturbation}) matches $\delta \approx 10^{-7}, 10^{-8}, 10^{-9}$, respectively.


\begin{table}[h!]
	\centering
	\caption{ The results of recovering derivative $F^{(2,0)}$ for noise from quadrature formula }
	\label{tbl3}
	\begin{tabular}{|c|c|c|c|}
		\hline
		$\delta$ & $10^{-7}$ & $10^{-8}$ & $10^{-9}$  \\ \hline
		Error L2        &  $4,5 \cdot 10^{-6}  $       &  $ 9,65 \cdot 10^{-7}   $  &  $8.6  \cdot 10^{-8}   $  \\ \hline
		Error C        &  $3,18 \cdot 10^{-5} $       &  $4,2 \cdot 10^{-6}   $  &  $4  \cdot 10^{-7}   $  \\ \hline
		n  &  19     & 31    & 43    \\ \hline
		h  &  $8\cdot 10^{-5}$    & $2\cdot 10^{-5} $  & $ 8\cdot10^{-6}$    \\ \hline
	\end{tabular}
\end{table}

The results of the numerical experiment are shown in Table \ref{tbl3} and Graph \ref{Fig3}.

\begin{figure}[h!]
	\begin{minipage}[h]{0.5\linewidth}
		\center{\includegraphics[width=1\linewidth]{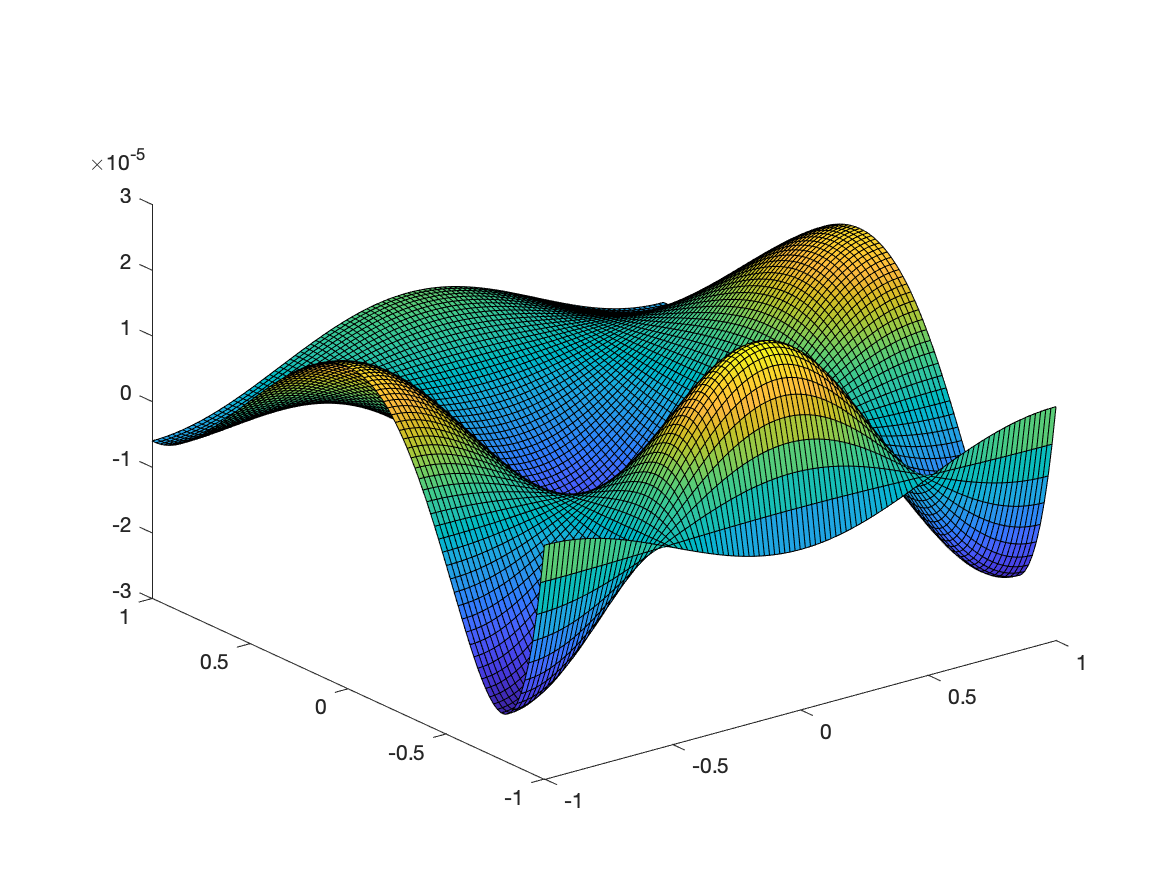} \\ a)}
	\end{minipage}
	\begin{minipage}[h]{0.5\linewidth}
		\center{\includegraphics[width=1\linewidth]{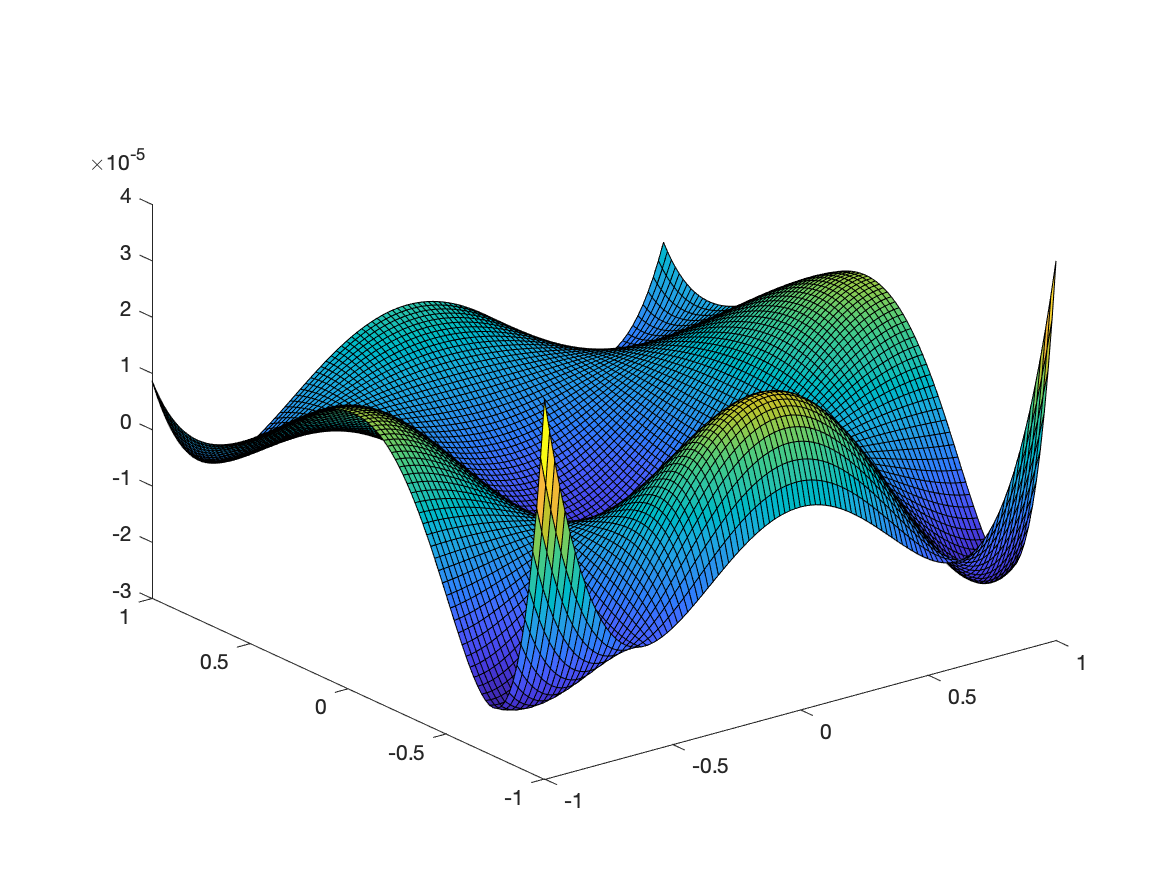} \\ b)}
	\end{minipage}
	\hfill
	\begin{minipage}[h]{0.5\linewidth}
		\center{\includegraphics[width=1\linewidth]{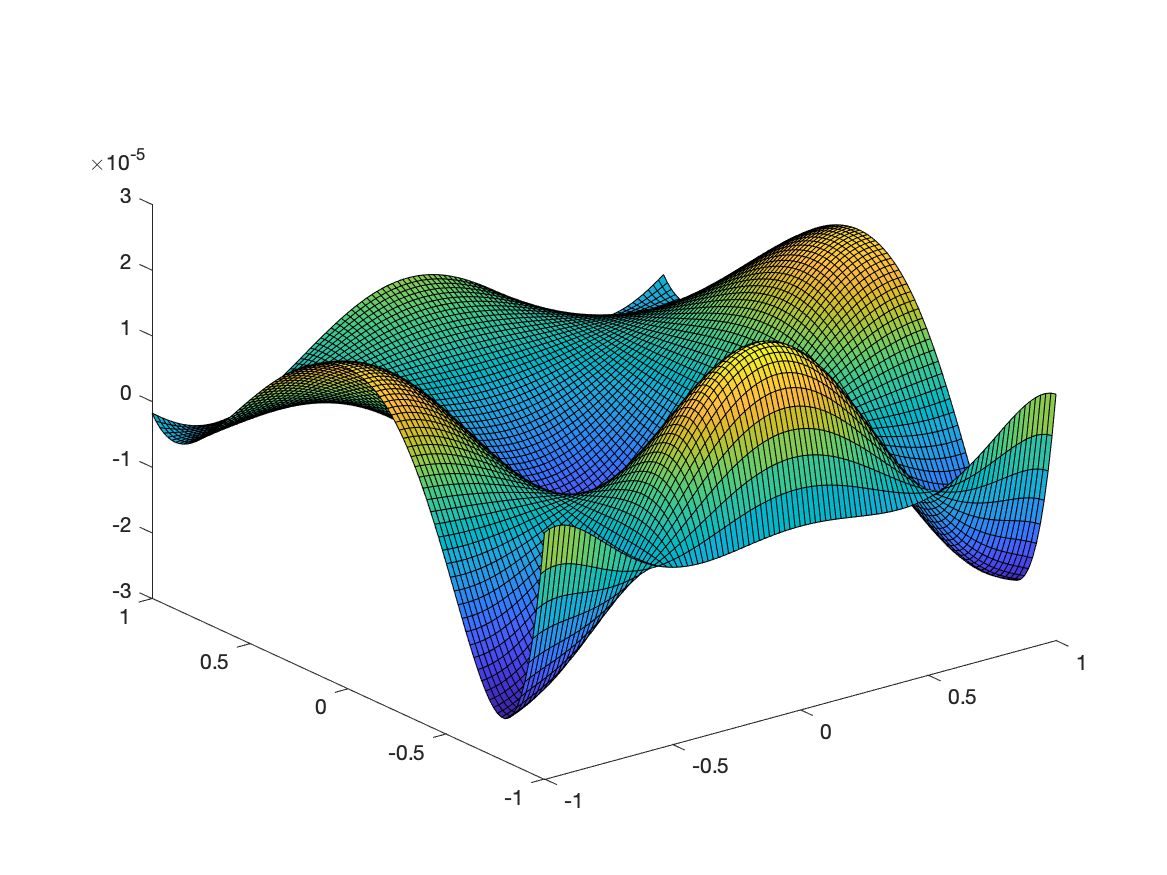} \\ c)}
	\end{minipage}
	\begin{minipage}[h]{0.5\linewidth}
		\center{\includegraphics[width=1\linewidth]{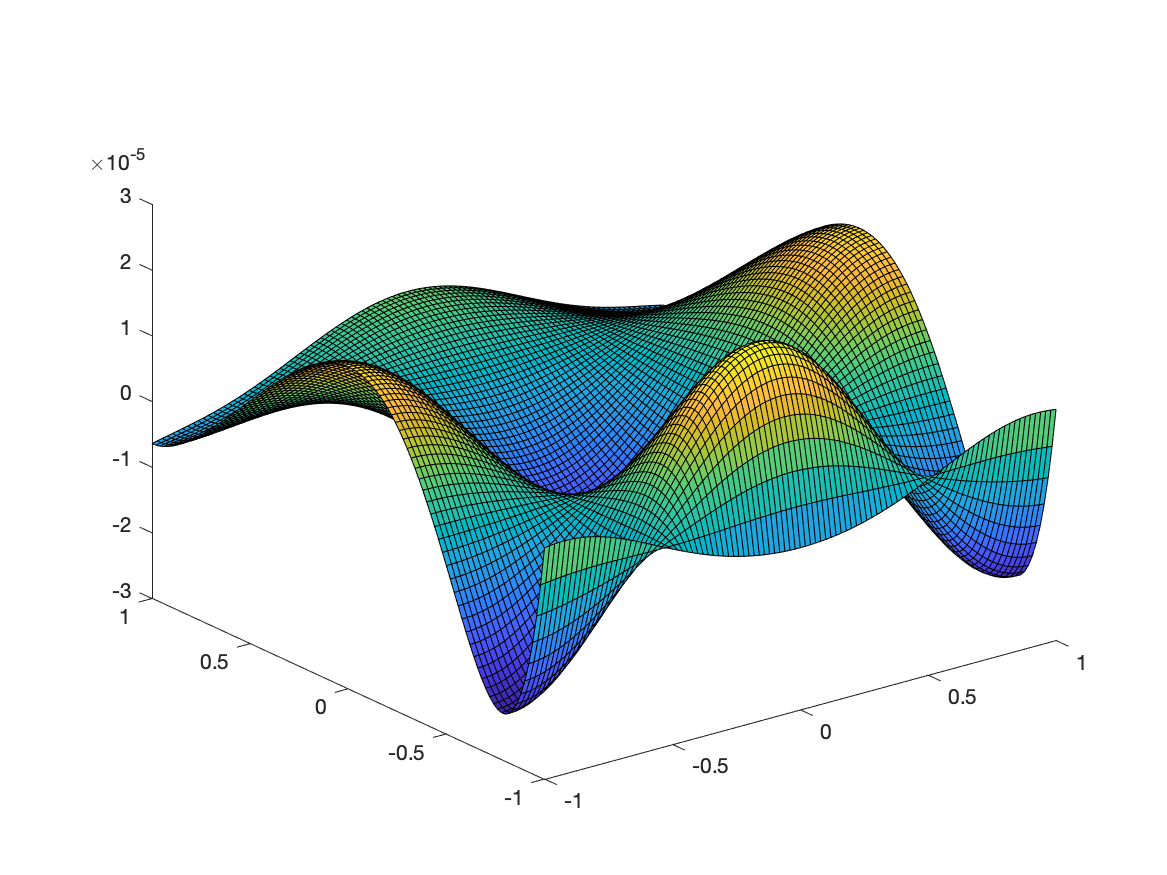} \\ d)}
	\end{minipage}
	\caption{Recovery of the derivative $F^{(2,0)}$ with noise in the input data arising from the quadrature formula. The exact derivative $F^{(2,0)}$ (Fig. a) ); approximation to $F^{(2,0)}$ for    $\delta= 10^{-7}$ (Fig.  b) );    for $\delta= 10^{-8}$ (Fig. c) ) and $\delta= 10^{-9}$ (Fig.  d) ). }
	\label{Fig3}
\end{figure}

 \section{ Acknowledgements }
This project has received funding through the MSCA4Ukraine project,
which is funded by the European Union.
In addition,
the first named author is supported by the Volkswagen Foundation project "From Modeling and Analysis to Approximation".
Also, the authors acknowledge partial financial support due to the project "Mathematical modelling of complex dynamical systems and processes caused by the state security" (Reg. No. 0123U100853).

\begin{small}
\begin{flushleft}
\textsc{Institute of Mathematics, National Academy of Sciences of Ukraine,
3, Tereschenkivska Str., 01024, Kiev, Ukraine}\\
\textit{E-mail address}: semenovaevgen@gmail.com, solodky@imath.kiev.ua
\end{flushleft}
\end{small}


\begin{thebibliography}{99}
\begin{small}
	
\bibitem{Ahn&Choi&Ramm_2006}
Ahn~S. \ A scheme for stable numerical differentiation / S.~Ahn, U.J.~Choi,  A.G.~Ramm
// J. Comput. Appl. Math. -- 2006.-- Vol.~
186(2).-- P. 325--334



\bibitem{Cul71} J. Cullum, Numerical Differentiation and Regularization.
SIAM Journal on Numerical Analysis. -- 1971. -- Vol.8. -- P. 254--265.

\bibitem{Groetsch_1991} Groetsch, C. W. Differentiation of Approximately Specified Functions. The American Mathematical Monthly. --
1991. -- Vol.98(9) -- P. 847--850.

\bibitem{Groetsch_1992_V74_N2}
 Groetsch~C.W.\
 Optimal order of accuracy in Vasin's method for differentiation of noisy functions
 /C.\,W.~Groetsch
 //J. Optim.Theory Appl.-- 1992.-- Vol. 74(2). -- P. 373--378.


\bibitem{HamKan2018} Hamarik, U., Kangro, U. On Self-regularization of Ill-Posed Problems in Banach Spaces by Projection Methods. In:Hofmann, B., Leitao, A., Zubelli, J. (eds) New Trends in Parameter Identification for Mathematical Models. 
    Trends in Mathematics. Birkhauser, Cham. --2018. -- P. 89--105.

\bibitem{Hanke&Scherzer_2001_V108_N6}
 Hanke~M.\
 Inverse problems light: numerical differentiation
/M.~Hanke, O.~Scherzer
//  Amer. Math. Monthly.-- 2001.-- Vol.~108(6).-- P. 512--521.

%

\bibitem{Meng&Zhaoa&Mei&Zhou_2020}
Meng ~Z. \ Numerical differentiation for two-dimensional functions by a Fourier extension method
/Z.~Meng, Z.~Zhaoa,
D.~Mei, Y.~Zhou
// Inverse Problems in Science and Engineering.-- 2020.-- Vol.~28(1).-- P. 1--18.

\bibitem{Lu&Naum&Per}
Lu ~S.\ Legendre polynomials as a recommended basis for numerical differentiation in the presence of stochastic white
noise / S. Lu, V. Naumova, S.\, V. Pereverzev
// J. Inverse Ill-Posed Probl.-- 2013.-- Vol.~21(2).-- P. 193--216.


\bibitem{Nakamura&Wang&Wang_2008}
Nakamura~G.\
 Numerical differentiation for the second order derivatives of functions of two
variables, /G. Nakamura, S.\, Z. Wang, Y.\, B. Wang
// J. Comput. Appl. Math.-- 2008.-- Vol.  212(2).-- P. 341--358.


\bibitem{Zhao_2010}
Zhao~Z.\
 A truncated Legendre spectral method for solving numerical differentiation
 /Z.~Zhao
 //  International Journal of
Computer Mathematics.-- 2010.-- Vol.~87.-- P. 3209--3217.

\bibitem{Zhao&Meng&Zhao&You&Xie_2016}
Zhao~Z.  \ A stabilized algorithm for multi-dimensional numerical differentiation /Z.~Zhao, Z.~Meng, L.~Zhao, L.~You,
O.~Xie
// Journal of Algorithms and Computational Technology.-- 2016.-- Vol.~10(2).-- P. 73--81.

\bibitem{Dolgopolova&Ivanov_USSR_Comput_Math_Math_Phys_1966_Eng}
Dolgopolova~T.F. \
 On numerical differentiation
 /T.F.~Dolgopolova, V.K.~Ivanov
 // Zh. Vychisl. Mat. and Mat. Ph.-- 1966. -- Vol.~6(3).-- P. 223--232.

\bibitem{EgorKond_1989}
Yu. V. Egorov, V.A. Kondrat'ev, On a problem of numerical differentiation, Vestnik Moskov. Univ. Ser. I Mat. Mekh. 3
(1989), P. 80-81.

\bibitem{Mileiko_Solodkii_2014}
Solodky~S.G.\ The minimal radius of Galerkin information for severely ill-posed problems /S.G.~Solodky, G.L.~Myleiko
//  Journal of
Inverse and Ill-Posed Problems.-- 2014.-- Vol.~22(5).-- P. 739--757.
%
\bibitem{Mileiko_Solodkii_2017_UMJ}
Mileyko~G.L. \
  Hyperbolic cross and complexity of different classes of linear ill-posed problems
 / G.L.~Mileyko, S.G.~Solodky
 // Ukr. Mat. J.-- 2017.-- Vol.~69(7).-- P. 951--963.

\bibitem{Mul69}
 M\" uller~C. \
 Foundations of the Mathematical Theory of Electromagnetic Waves
 /C.~ M\" uller --
 Springer--Verlag, Berlin,
Heidelberg, New York, 1969.


\bibitem{Pereverzev_Computing_1995} 
 Pereverzev~SV. Optimization of projection methods for solving ill-posed problems. Computing. -- 1995.
 -- Vol.55. P. 113--124.

\bibitem{Sem_Sol_2008}
Lebedeva E.V., Solodky S.G. \ Approximation of finite-section equations by piecewise constant functions //
Computational Mathematics and Mathematical Physics. -- 2008. -- Vol.48(5). -- P.693-706.

\bibitem{Sem_Sol_2021}
Semenova Y.V. \ Error bounds for Fourier-Legendre truncation method in numerical differentiation, /Semenova Y.V.,
Solodky S.G.
// Journal of Numerical and Applied Mathematics.-- 2021.-- Vol.~137(3).--P.113--130.

\bibitem{Sol_Stas_JC2020} S.\,G. Solodky, S.\, A. Stasyuk, Estimates of efficiency for two methods of stable numerical summation of
smooth functions, {\it Journal of Complexity} {\bf 56} (2020).  https://doi.org/10.1016/j.jco.2019.101422



\bibitem{Sol_Stas_UMZ2022}
Solodky S.G. On optimization of methods of numerical differentiation for bivariate functions, / Solodky S.G., Stasyuk S.
// Ukr. Mat. J.-- 2022.-- Vol.~74(2).-- P. 253--273.


\bibitem{Qian&Fu&Xiong&Wei_2006}
Qian~Z.\
 Fourier truncation method for high order numerical derivatives
 /Z.~Qian,  C.L.~Fu,  X.T.~Xiong,  T.~Wei,
 // Appl.
Math. Comput.-- 2006.-- Vol.~181(2).-- P. 940--948.

\bibitem{PS1996}
S.\, V. Pereverzev, S.\, G. Solodky, The minimal radius of Galerkin information for the Fredholm problem of the first
kind, {\it Journal of Complexity} {\bf 12}~(4) (1996), P. 401--415.

\bibitem{Ramm_1968_No11}
 Ramm ~ A.G., \
 On numerical differentiation
 / A. \, G. ~ Ramm
 // Izv. Vuzov. Matem.-- 1968.-- Vol.~11.-- P. 131--134.

\bibitem{RammSmir_2001}
A.G. Ramm, A.B. Smirnova On stable numerical differentiation Math. Comput., 70 (2001), P. 1131-1153

\bibitem{SSS_CMAM}
Semenova Y.V.\ Application of Fourier truncation method to numerical differentiation for two variables functions
/Y.\,V. Semenova, S.\, G. Solodky, S.\, A. Stasyuk
// Computational Methods in Applied Mathematics, vol. 22, no. 2, 2022, pp.
477-491. https://doi.org/10.1515/cmam-2020-0138

\bibitem{SolSem2012} {\it  Solodky S.G., Semenova E.V. } On optimal order accuracy of solving  Symm's integral equation // Zh. Vychisl. Mat. i Mat. Fiz. 2012. V. 52, No 3. P. 472-482



\bibitem{TrWW} J.\,F. Traub, G.\,W. Wasilkowski, H. Wozniakowski, {\it Information-Based Complexity}, Academic
Press, New York, 1988.

\bibitem{TrW} J.\,F. Traub,  H. Wozniakowski, A General Theory of Optimal Algorithms, Academic Press, New York,
1980.

\bibitem{Wang_Hon_Ch_2006} Y.B. Wang, Y.C. Hon, J. Cheng,
Reconstruction of high order derivatives from input data // J. Inverse Ill-posed Probl., 14 (2006), P. 205-218.


 \bibitem{VainHam} {\it  Vainikko, G. M.;  Khyamarik, U. A.}  Projection methods and self-regularization in ill-posed problems.  // (Russian) Izv. Vyssh. Uchebn. Zaved. Mat. 1985. V. 29. P. 1-17.

\bibitem{VasinVV_1969_V7_N2}
Vasin~V.V., \
 Regularization of the numerical differentiation problem
 /V.V. ~Vasin
 // Mat. app. Ural un-t.-- 1969.-- Vol.~7(2).-- P. 29--33.



%

%

\end{small}
\end{thebibliography}
\end{document}